\documentclass[reprint,superscriptaddress,showpacs,amsmath,amssymb,aps,prx,floatfix,]{revtex4-2}

\usepackage[utf8]{inputenc}	
\usepackage{amsmath}
\usepackage{graphicx}		
\usepackage{natbib}
\usepackage{textcomp}
\usepackage{gensymb}
\usepackage[hidelinks]{hyperref}
\usepackage{xcolor}
\usepackage{siunitx}
\usepackage{mathrsfs}
\usepackage{multirow}
\usepackage{xfrac}
\usepackage{float}
\usepackage{mathtools}
\usepackage{comment}
\newcommand*{\meanR}{\ensuremath{R}}
\newcommand*{\walrus}{\ensuremath{\coloneqq}}
\newcommand{\kl}[1]{\textcolor{black}{#1}}

\begin{document}

\title{Small changes at single nodes can shift global network dynamics}

\author{Kalel L. Rossi}
\affiliation{Theoretical Physics/Complex Systems, ICBM, Carl von Ossietzky University of Oldenburg, Oldenburg, Lower Saxony, Germany}
\author{Roberto C. Budzinski}
\affiliation{Department of Mathematics, Western University, London, Ontario, Canada}
\affiliation{Brain and Mind Institute, Western University, London, Ontario, Canada}
\affiliation{Western Academy for Advanced Research, Western University, London, Ontario, Canada}
\author{Bruno R. R. Boaretto}
\affiliation{Institute of Science and Technology, Federal University of São Paulo, São José dos Campos, São Paulo, Brazil}
\author{Lyle E. Muller}
\affiliation{Department of Mathematics, Western University, London, Ontario, Canada}
\affiliation{Brain and Mind Institute, Western University, London, Ontario, Canada}
\affiliation{Western Academy for Advanced Research, Western University, London, Ontario, Canada}
\author{Ulrike Feudel}
\affiliation{Theoretical Physics/Complex Systems, ICBM, Carl von Ossietzky University of Oldenburg, Oldenburg, Lower Saxony, Germany}

\begin{abstract}
Understanding the sensitivity of a system's behavior with respect to parameter changes is essential for many applications. This sensitivity may be desired - for instance in the brain, where a large repertoire of different dynamics, particularly different synchronization patterns, is crucial - or may be undesired - for instance in power grids, where disruptions to synchronization may lead to blackouts. In this work, we show that the dynamics of networks of phase oscillators can acquire a very large and complex sensitivity to changes made in either their units' parameters or in their connections - even modifications made to a parameter of a single unit can radically alter the global dynamics of the network in an unpredictable manner. As a consequence, each modification leads to a different path to phase synchronization manifested as large fluctuations along that path. This {\it dynamical malleability} occurs over a wide parameter region, around the network's two transitions to phase synchronization. One transition is induced by increasing the coupling strength between the units, and another is induced by increasing the prevalence of long-range connections. Specifically, we study Kuramoto phase oscillators connected under either Watts-Strogatz or distance-dependent topologies to analyze the statistical properties of the fluctuations along the paths to phase synchrony. We argue that this increase in the dynamical malleability is a general phenomenon, as suggested by both previous studies and the theory of phase transitions.
\end{abstract}

\maketitle

\section{Introduction}\label{sec:intro}
Several systems of practical and theoretical importance are composed of, or can be modeled as, networks of interacting units. Examples from different research areas include power grids (networks of producers and consumers of electrical energy) \cite{motter2013spontaneous}, food webs \cite{dunne2002foodweb}, networks of electronic elements \cite{crotty2010josephson}, coupled lasers \cite{nixon2011synchronized} and neurons in the brain \cite{varela2001brainweb}. An important question is how the dynamics of the single units in the network impact the overall dynamics of the system, and what happens if these units are modified. This could be by changing the units' parameters - e.g. in ecological systems, what happens if the reproduction rate of a prey increases? in power grids, can a change in the parameters of a single generator cause a large disruption, such as a blackout? or the modification could be by shocking the units into a different dynamical state - e.g. in the brain, how can an epileptic seizure be stopped by employing a current pulse in one particular brain region? \kl{A regime in which changes in a single unit can alter the whole network's behavior can thus be either dangerous or advantageous, and is an important topic of research, which we address in this work.}

In both power grids and the brain, an important phenomenon is synchronization, i.e. the coherence of frequencies or even phases of oscillations: it is e.g. crucial for power grids to have their elements synchronized in the \SI{50/60}{\hertz} regime \cite{witthaut2022collective}. Moreover, several functional roles have been ascribed to synchronization in the brain \cite{fries2015rhythms, varela2001brainweb, singer1999neuronal}.  Particularly for systems in which synchronization is an essential process for functioning, the question of sensitivity with respect to perturbations becomes important. This has been recognized in the literature, and various types of perturbations - to \kl{the system's dynamical state}, to topology, and to parameters of the units - have been considered to study the vulnerability either of the synchronized state itself or of the transition to synchronization \cite{pikovsky2001synchronization, arenas2008synchronization}. Especially interesting are studies devoted to the impact of perturbations to only a small part of the network, to even a single node or a single link. For networks of identical units it has been shown that perturbations to the state of a single node can lead to desynchronization of parts or even the whole network \cite{mitra2017multiple, medeiros2018boundaries}. The same outcome can be observed by adding/removing a single link in the network \cite{witthaut2012braess, coletta2016linear, mihara2022sparsity}. Similarly, studies in networks of non-identical units have proposed measures and identified topological properties responsible for increased vulnerability in networks due to perturbations in the state \cite{menck2014how, halekotte2020minimal} or parameters \cite{manik2017network} of the units.

\begin{figure}[htb]
    \centering
    \includegraphics[width=0.9\columnwidth]{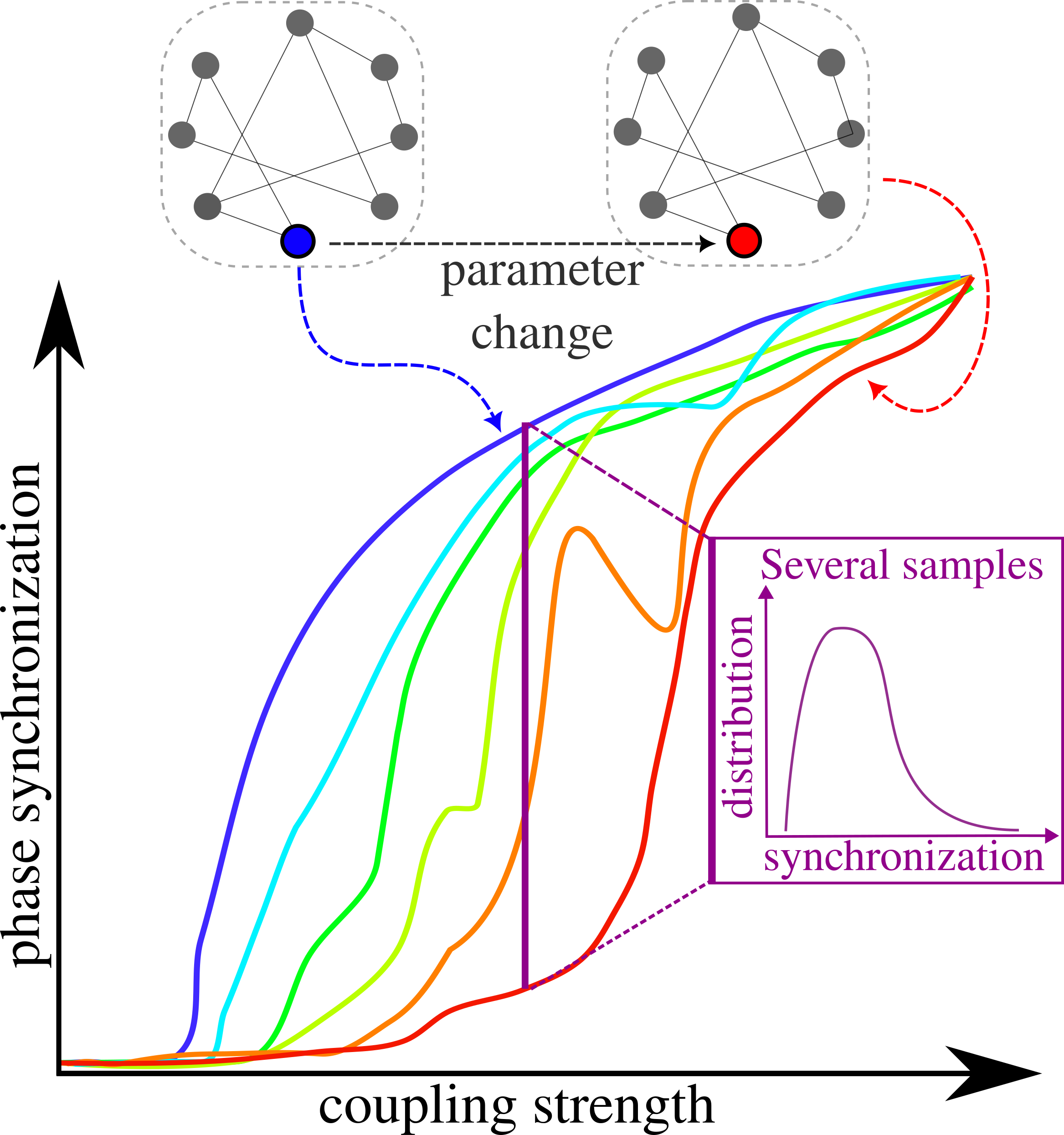}
    \caption{Illustration of sample-to-sample fluctuations (dynamical malleability) in a classical transition to phase synchronization. Each path to phase synchronization corresponds to a different sample, which can differ from the others because of the difference in the parameter of a single unit. We see the transitions are different, as the critical value of coupling strength and the profile of the transition differ. Fixing the coupling strength, we can also look at the distribution of the degree of phase synchronization across samples (purple inset). The magnitude of the sample-to-sample fluctuations peaks during the transitions, at the edge of phase synchronization \cite{buendia2022the}.}
    \label{fig:sketchmalleability}
\end{figure}

In this work we focus on perturbations to parameters controlling the dynamics of the units. The networks we study have two distinct transitions to phase synchronization. \kl{In both, there is a large increase in the system's sensitivity, such that changes to parameters of even single units can radically alter the dynamics of the whole network, including its transition to phase synchronization and its spatiotemporal patterns} (cf. illustrated in Fig. \ref{fig:sketchmalleability}). This \textit{dynamical malleability} can cause problems in real systems in two major ways: (i) the large magnitude of the fluctuations in the dynamics, which can lead to drastic changes in behavior, and which makes the exact path to phase synchronization hardly predictable and (ii) the complexity of the fluctuations: it is not clear in principle which units have to be chosen, or how large the change in their parameters must be, to cause the largest deviation from a mean path to phase synchronization; also conversely, it is not clear which chosen units or changes in parameters can keep the network in a similar synchronization state. This clearly important issue for the design and control of systems motivates our study to analyze the mechanisms leading to those large fluctuations. 

\kl{To address this issue, we study networks of Kuramoto phase oscillators organized in a ring lattice. They constitute a paradigmatic model for synchronization \cite{kuramoto1975self,acebron2005kuramoto,rodrigues2016the} and have been established as a model for} real-world systems like the brain \cite{ponce-alvarez2015restingstate, cabral2011role, rodrigues2016the}, Josephson junctions \cite{josephson1964coupled, crotty2010josephson} and chemical oscillators \cite{marek1975synchronization, neu1979chemical}. These phase oscillators are limit cycle oscillators for which it is assumed that the amplitude dynamics can be neglected due to a strong convergence compared to the coupling strength and they can, hence, be described only by a phase (an angle) that tends to evolve linearly according to a natural frequency. In the network, the coupling between two units is realized by the sinusoidal difference of their phases. Part of the appeal for this model is that the equations are generic, in the sense that they describe the behavior of any system of limit-cycle, slightly non-identical, oscillators with weak coupling \cite{kuramoto1984chemical, strogatz2000from, kuramoto2019on}. They are connected here in either of two classes of network topologies: one obtained by a successive random rewiring of connections starting from a $k$-nearest-neighbor ring (Watts-Strogatz, small-world networks \cite{watts1998collective}) up to a random network; and another in which all units are connected, but the connection weights decay with the distance between nodes according to a power law \cite{rogers1996phasetransitions}. The latter includes the globally coupled mean-field topology in the limit of no decay. The two classes have similarities: firstly, they have similar $k$-nearest-neighbor rings; and secondly, they include networks of mean-field type (random or globally coupled) \cite{skardal2020higher, hong2013link}. They also have differences: the first class is sparsely connected, the other densely; the first has link-disorder (different rewirings lead to different networks), the second does not. Despite these differences, however, both topologies present a similar phenomenology: going from $k$-nearest-neighbor rings to mean-field rings leads to a transition from desynchronization to phase synchronization \cite{hong2002synchronization, rogers1996phasetransitions}. Furthermore, the systems become very sensitive to parameters changes (very dynamically malleable) in the path to phase synchronization (at the edge of phase synchronization \cite{buendia2022the}). The parameter changes considered here are mainly changes to the oscillators' natural frequencies, but we also show that changes to the topology have similar effects.

This large dynamical malleability can be connected to statistical physics, as transitions to phase synchronization have been well-established as non-equilibrium phase transitions, especially for systems like Kuramoto oscillators under global coupling \cite{kuramoto1984chemical,hong2007entrainment, peter2018transition}. However, the comparison to phase transitions is rigorously true only in the limit of an infinite number of oscillators (thermodynamic limit), such as the one originally studied by Kuramoto \cite{kuramoto1984chemical}, since phase transitions can rigorously only occur in this limit \cite{brankov2000theory}.  For finite systems, such as the ones we study here and any system that occurs in practice, a behavior analogous to a phase transition occurs, in which the system's behavior changes as a control parameter (e.g. coupling strength) is changed, but does so over a wider parameter interval, instead of a single point: the transition becomes blurred \cite{brankov2000theory, binder1987finite}. For finite systems with quenched disorder - meaning heterogeneous and constant parameters of the units, such as the frequency of the oscillators - each particular realization of the parameters (sample) also transitions at different intervals of the control parameter. This results in large \textit{sample-to-sample (STS) fluctuations} (i.e. dynamical malleability) near the phase transition \cite{sornette2006critical, wiseman1995lack}. The increase of sample-to-sample fluctuations near phase transitions has been described for Kuramoto oscillators under a few topologies \cite{peter2018transition, hong2007entrainment, hong2007finitesizescalingpre}. For globally coupled Kuramoto oscillators, Hong et al. \cite{hong2007entrainment} offer a statistical analysis including the STS fluctuations, which increase significantly during the transition to phase synchronization. For random connectivity, Hong et al. \cite{hong2013link} offer a similar treatment considering disorder in the connections. These help consolidate the generality of sample-to-sample fluctuations, which can be expected whenever a system has quenched disorder and is near a phase transition.

Despite the commonness of the behavior, its potentially large effects on the system's dynamics are still poorly understood. In this work, we highlight two specific points: the first is the large magnitude of the STS fluctuations, and the corresponding \textit{large sensitivity} that emerges from them - as we mentioned, changing the frequency of a single unit in the network may radically alter the whole system's behavior. The second point is the \textit{complexity} of the changes. Several methods have been developed to investigate transitions to phase synchronization in networks \cite{peter2018transition, skardal2014optimal, brede2008synchrony, carareto2009optimized}. We have tested them, but none could satisfactorily describe the STS fluctuations we observe. 
Firstly, there is a mechanism for STS fluctuations, proposed by Peter and Pikovsky \cite{peter2018transition}, in which different realizations of the natural frequencies, with fixed mean and standard deviation, can lead to (i) most units close to the mean, with a few extreme outliers; or to (ii) some units spread around the mean, with only moderate outliers. This is quantified by the distribution's kurtosis. The transition for synchronization due to increasing coupling strength starts earlier for the former (higher kurtosis), as most units have similar frequencies, but finishes later, as the extreme outliers need a higher coupling strength to synchronize. This mechanism, however, is not capable of explaining the STS fluctuations we observe here, as we discuss in the conclusions. Second, the synchrony alignment function, shown analytically by \cite{skardal2014optimal} to be related to the degree of phase synchronization in the limit of strong synchrony, does not work for the whole range of phase synchronization we observe here, and thus does not explain the whole behavior. Thirdly, other measures that have been observed in the literature to correlate to phase synchronization do not work in the malleable networks. These are: (i) the proportion $p_-$ of links connecting nodes with natural frequencies of different signs \cite{brede2008synchrony};  (ii) the correlation $c_\omega$ between the oscillators' natural frequencies, taking into account the connectivity of the network \cite{brede2008synchrony, carareto2009optimized}; (iii) the correlation between natural frequencies and the node's number of connections (degree) \cite{skardal2014optimal}; and (iv) the correlation between the average frequency between neighbors of a node and the node's own frequency \cite{skardal2014optimal}.

As illustrated in Fig. \ref{fig:sketchmalleability}, the degrees of phase synchronization of the samples, for fixed coupling strength, have a certain distribution. While works on statistical physics of Kuramoto networks generally assume that these STS fluctuations are normally distributed \cite{tang2011finite, hong2007entrainment, hong2015finite}, evoking the central limit theorem, we show here that for our networks, up to even size $N=5000$, the distributions depend considerably on the coupling strength and the topology, and are usually far from being normal, even in globally coupled networks. 

\kl{Besides changing the units' parameters, we also change their initial conditions, which gives insight into the system's response to perturbations altering its dynamical state. We show that the number of attractors increases during the transition to phase synchronization in the Watts-Strogatz networks. This, we believe, is a novel result in the literature. The increased multistability acts as a dynamical mechanism to increase the dynamical malleability, but is not necessary, as the distance-dependent networks seem to be monostable.} This further highlights the importance of topology, which we quantify via the ratio between short-range and long-range connections in the networks. We demonstrate with this that the dynamical malleability - and the multistability, for WS networks - peaks over specific ranges of this ratio, for a certain balance between short and long-range connections. 

With this work, we therefore hope to demonstrate the importance of dynamical malleability and the large sensitivity it brings, and to encourage further theoretical advancements in this area, which are needed to properly describe the wide range of behaviors and to offer tools for practical applications.  

\section{Methodology}\label{sec:methodology}

In the Kuramoto model \cite{kuramoto1975self, kuramoto1984chemical}, each oscillator is described by a phase which evolves in time according to
\begin{equation}
    \Dot{\theta_{i}} = \omega_{i} + \epsilon \sum\limits_{j=1}^{N} A_{i,j} \sin{(\theta_{j} - \theta_{i})},
    \label{eq:main_kuramoto}
\end{equation}
where $\theta_{i}(t)$ is the phase of the $i$-th oscillator at time $t$, $\omega_{i}$ is its natural frequency, $\epsilon$ is the coupling strength, $N$ is the number of oscillators, and $A_{i,j}$ is the $(i,j)$-th element of the adjacency matrix $\bm{A}$. Throughout this work, we initially draw each frequency randomly from a Gaussian distribution with mean $\mu = 0.0$ and standard deviation $\sigma = 1.0$, generating a sequence $\{\omega_{i}\}, i=1,\cdots,N$. Then, different realizations can (i) shuffle these frequencies, generating another sequence $\{\omega_{i}\}_{\mathrm{shuffled}} = \mathrm{shuffle}(\{\omega_{i}\})$; or (ii) switch one selected unit's frequency to another value $\omega_{\mathrm{new}}$. 

The networks in this work are coupled in a ring lattice of $N=501$ units with periodic boundary conditions, and follow one of two types of topology. 
The first is the Watts-Strogatz topology \cite{watts1998collective}, which interpolates between a regular and a random topology with a parameter $p$, the rewiring probability: at one extreme ($p = 0$), the topology is a $k$-nearest-neighbor lattice; starting from it, connections are randomly chosen, according to the probability $p$, and rewired to another randomly chosen connection. Doing this, the networks have a significant decrease in the mean distance between nodes, but remain very clustered, generating small-world topologies. The other extreme ($p = 1$) is then a random topology. These networks are unweighted so their adjacency matrix's elements are $A_{i,j} = 1$ if $i$ and $j$ are connected, and $0$ otherwise. 
The second class of networks follow a distance-dependent powerlaw scheme, in which any given node receives connections with weights decaying based on the distance to that node. Each element of the adjacency matrix is $A_{ij} =\frac{1}{\eta(\alpha) (d_{ij})^\alpha}$, where $d_{ij}$ is the edge distance between oscillators $i$ and $j$, defined as  $d_{ij} = \mathrm{min}(|i - j|, N - |i - j|)$, and $\eta(\alpha)$ is a normalization term given by: $\eta(\alpha) = \sum\limits_{j=1}^{N^\prime}\dfrac{2}{j^{\alpha}}$, such that the temporal evolution of the phases can be written as:
\begin{equation}
    \Dot{\theta_{i}} = \omega_{i} + \frac{\epsilon}{\eta(\alpha)} \sum\limits_{j=1}^{N^\prime} \frac{1}{j^\alpha} \left[ \sin{(\theta_{i+j} - \theta_{i})} + \sin{(\theta_{i-j} - \theta_{i})} \right],
    \label{eq:main_kuramoto_powerlaw}
\end{equation}
where $N^\prime = \frac{N-1}{2}$ denotes half the amount of units that $i$ is connected to (one half of the ring's length, discounting the unit $i$ itself). The equation explores the symmetry in the network to switch the summation across the network to a summation across only half, multiplied by $2$.
The powerlaw decay is thus controlled by $\alpha$, the locality parameter. For $\alpha = 0$, the network is globally coupled with equal weights between every node. As $\alpha$ increases, the weights are redistributed, so that closer units (in terms of edge-distance) have bigger weights. At the extreme of $\alpha \to \infty$, only first-neighbors are connected.

Integration was performed using the Tsitouras 5/4 Runge-Kutta (Tsit5) method for Watts-Strogatz networks, and an adaptive order adaptive time Adams Moulton (VCABM) method for distance-dependent networks. The integrator method was chosen for distance-dependent networks for increased simulation speed, and results were robust to different integration schemes. All methods used the DifferentialEquations.jl package \cite{rackauckas2016differential}, written in the Julia language \cite{bezanson2017julia}. Additional computational packages used were PyPlot \cite{hunter2007matplotlib} for plotting and DrWatson.jl \cite{datseris2020drwatson} for code management. The code used for simulations is accessible in the repository \cite{rossi2022github}, with the parameters used in the simulations. In particular, the control parameters we used ($\alpha$, $p$ and $\epsilon$) were generated to be equally separated inside their particular range (in the case of $p$, equally separated in a log scale). The values shown on all figures were thus not chosen by hand, and we have verified they correspond to typical behaviors of the systems.

The degree of phase synchronization of the network can be quantified by the circular average of the phases 
\begin{equation}
    r(t) = \frac{1}{N}\left|\sum\limits_{j=1}^N \exp{(i \theta_j(t))}\right|,
\end{equation} 
with $i = \sqrt{-1}$. The quantifier ranges from $0$ to $1$: if $r(t) = 1$, all the phases are the same, and the system is completely globally phase-synchronized; if $r(t) = 0$, each oscillator has a pair that is completely out-of-phase, and the system can be completely globally phase-desynchronized or in a twisted state. We typically describe networks by the temporal average $R \walrus \frac{1}{T} \sum_t r(t)$  of their phase synchronization, with $T$ being the total simulation time excluding transients.

\section{Results}\label{sec:results}
\subsection{Dynamical malleability}
\begin{figure*}[htb]
    \centering
    \includegraphics[width=0.95\textwidth]{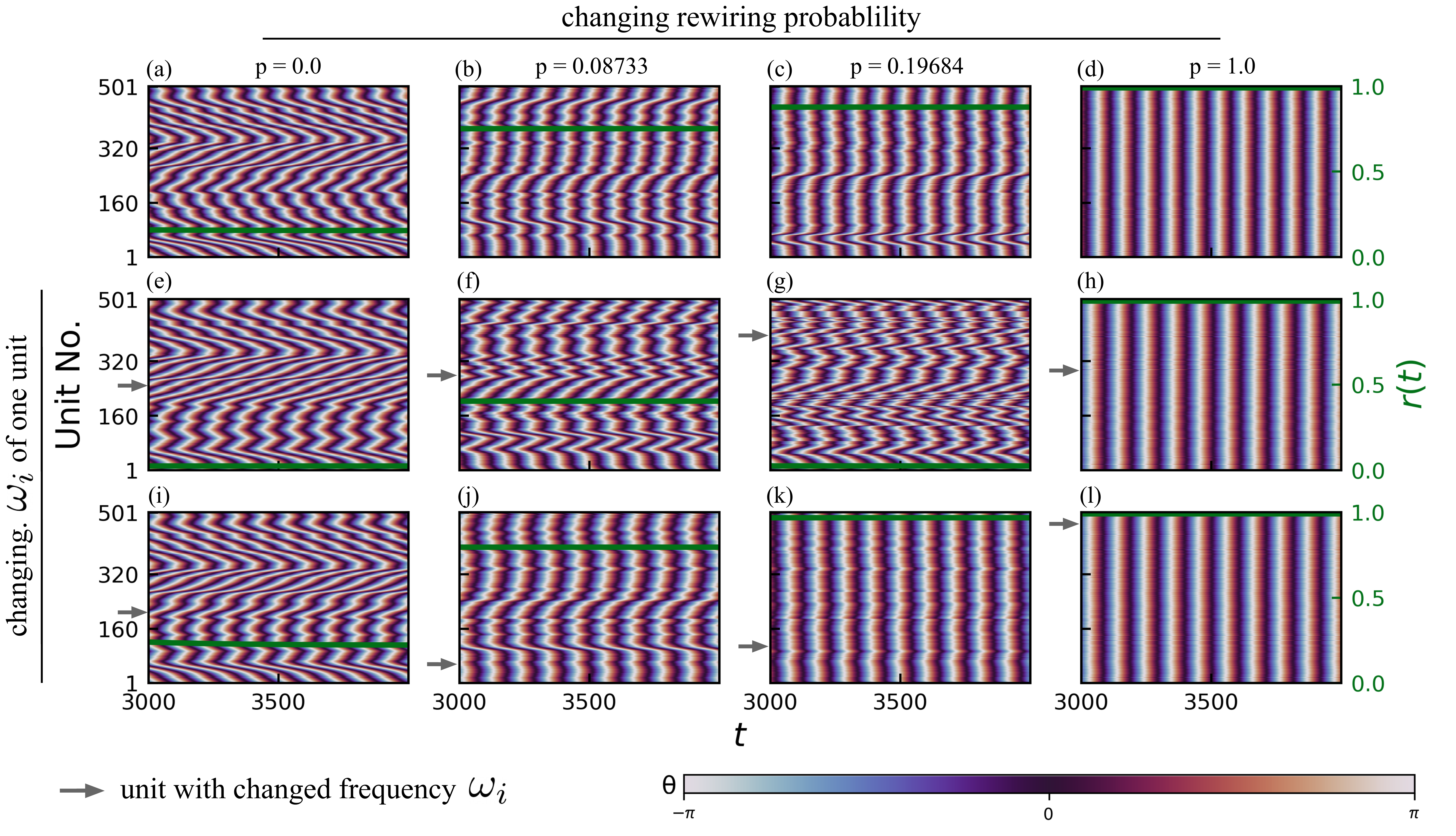}
    \caption{\textbf{Transition to phase synchronization and the effect of a single-unit change}. The figure shows the color-coded phases $\theta$ of all oscillators in the network and the degree of phase synchronization $r(t)$ (green line) across time for Watts-Strogatz networks. The coupling strength $\epsilon$ is fixed at $\epsilon = 4.51282$ and the natural frequencies ${\omega_i}$ in the first row are the same, generated by randomly drawing from a Gaussian distribution with zero mean and unitary standard deviation. Networks in the left column are two-nearest-neighbor lattices (rewiring probability $p = 0$); the short-range connections in these networks are then rewired in the following columns, with probability $p = 0.08733$ in the second column, $p = 0.19684$ in the third, and $p = 1.0$ in the fourth (leading to random networks). Increasing the proportion of long-range connections thus leads to more phase-synchronized networks.
    In the second and third rows, the natural frequency $\omega_i$ of a single unit $i$ (indicated by the gray arrows) is changed to a new value $\omega_i \to \omega_\mathrm{new} = 3.0$, with all other parameters being kept fixed, in particular the coupling strength $\epsilon$. The units shown in the figure were those which led to the smallest (second row) or highest (third row) degree of phase synchronization $\meanR$ out of all $N = 501$ units in the network for each value of $p$. Initial conditions were the same for all simulations, and were randomly drawn between $0$ and $2\pi$.
    }
    \label{fig:spatiotemporal}
\end{figure*}

The networks we study here, described in Eq. (\ref{eq:main_kuramoto}), have two characteristics: the oscillator's natural frequencies $\omega$, distinct for each unit, and their topology. To start this work, we introduce the dynamics of the network, and how they can change significantly if the natural frequencies of the units are slightly altered. 
For very small coupling strengths $\epsilon$ the oscillators are effectively uncoupled, and the phases oscillate without any significant correlation. As this $\epsilon$ increases, the instantaneous frequencies $\dot{\theta}_i$ align first, and the units' phases become locked, but not aligned (frequency but not phase synchronized). Whether the phases become aligned or not then depends on the topology \cite{medvedev2014small, hong2002synchronization}. In a two-nearest neighbor lattice, where only four nearby units are connected (two on each side), there is a topological limitation in the spread of interactions across the network which makes the oscillators arrange themselves in shorter-range patterns (Fig. \ref{fig:spatiotemporal}(a)) (an exception might occur if the coupling strength is extremely high, much bigger than the relevant values studied here).
As the short-range connections are randomly rewired to long-range connections, the shorter-range patterns give way to longer-range patterns, and the oscillators start to phase synchronize ((b) and (c)), until eventually a strong (though not complete) phase synchronization (PS) is reached (d). This occurs at different stages for each realization: for instance, panels (c) and (k) reach a high degree of PS, with the longer-range patterns, before panel (g) does. 
These changes occur due to the rewiring of connections, controlled by the random rewiring probability $p$ in the Watts-Strogatz (WS) networks.

The natural frequencies $\{\omega_i\} (i=1,\cdots, N)$ were kept constant across panels (a)-(d). Changing the frequencies, keeping initial conditions $\{\theta_i(0)\} (i=1,\cdots, N)$ fixed, leads to a different realization (also called sample), with possibly different dynamics. If the frequency of a single unit $\omega_i$ is changed to an arbitrary new value, for instance $\omega_\mathrm{new} = 3$, the network's behavior can be significantly altered (panels (e)-(h)). This is especially the case for networks with intermediate rewiring probabilities $p$, in which this single unit frequency change can bring the network from high to very low phase synchronization (panels (c) to (g)). The instantaneous frequencies typically remain synchronized, though their values might change. For random networks, phase synchronization is always maintained, though the instantaneous frequency values may change. 

Figure \ref{fig:spatiotemporal} thus illustrates that the long-term dynamics and phase synchronization differ in each realization. The realizations, created by changing the natural frequency of one unit or by changing the topology, are distinct dynamical systems, so it is not surprising to observe distinct long-term dynamics. It is, however, interesting to observe how these distinctions occur and depend on the topology. For instance, in networks of intermediate $p$, the distinctions occur such that the phase synchronization changes drastically; in random networks, they preserve the phase synchronization but alter the instantaneous frequencies. We also note that the behavior we describe is typical of the systems, and the values of $p$ and $\epsilon$ used here were generated as described in Section \ref{sec:methodology}. \kl{Since the fluctuations in the phase patterns (reflected in the phase synchronization) are clearer and more pronounced than the instantaneous frequency patterns, we now focus on the phase synchronization of the networks.}

\begin{figure*}[htb]
    \centering
    \includegraphics[width=0.95\textwidth]{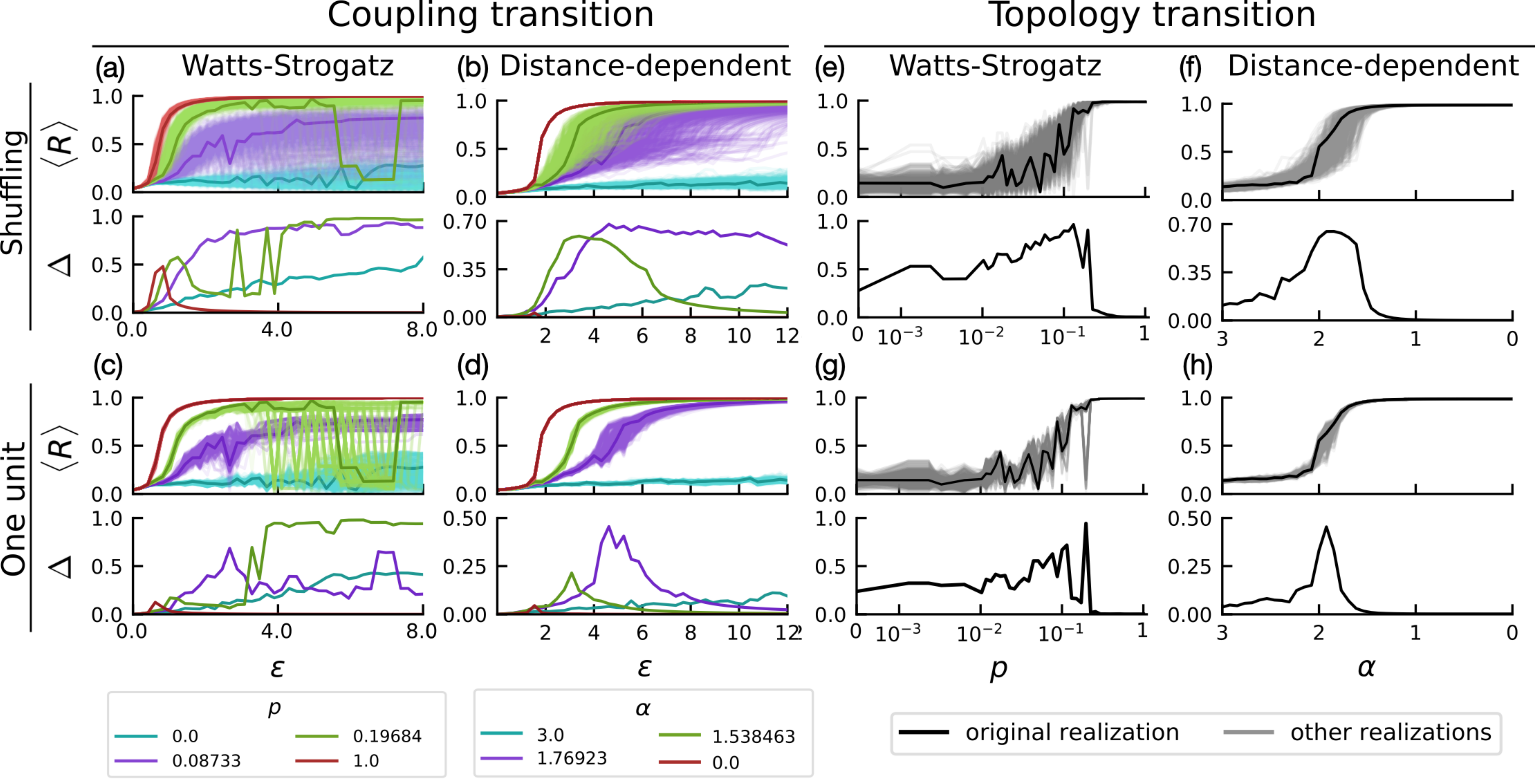}
    \caption{\textbf{Transitions to phase synchronization and sample-to-sample fluctuations}. Networks under Watts-Strogatz (WS) and distance-dependent (DD) topologies reach phase synchronization through either an increase in coupling strength $\epsilon$ (given the topology has a sufficient amount of long-range connections) or by switching short-range connections to long-range. Fluctuations in the degree of phase synchronization $\meanR$ between samples increase during the transitions, as can be seen by the differences in the same-colored curves and by $\Delta \coloneqq \meanR_\mathrm{max} - \meanR_\mathrm{min}$. Starting from a natural frequency sequence originally drawn from a Gaussian distribution (thicker lines), the other samples (thinner lines) can be generated by shuffling the natural frequencies or by switching the natural frequency of one unit to $\omega_\mathrm{new} = 3$. For intermediate networks (purple and green curves), the increase in fluctuations extends for a wide range of parameters and becomes considerably large. Each panel contains $501 = N$ realizations, with rewiring probabilities fixed for the coupling transition, with values shown in the legend, and coupling strength fixed in the topology transition at $\epsilon = 4.51282$ for WS and $\epsilon = 6.46154$ for DD. The initial conditions are the same across all realizations, and are randomly distributed from $0$ to $2\pi$. The curves of $\Delta$ are qualitatively similar with other dispersion measures, such as standard deviation, a possible difference being that the curves may be slightly shifted, as the measures can peak at slightly different values of the control parameter. We remark that the parameter values used in the simulations were generated as described in Section \ref{sec:methodology}, and were verified to correspond to typical parameter values.}
    \label{fig:transition_sync}
\end{figure*}

To obtain a comprehensive picture we now study an ensemble of samples obtained by shuffling the frequencies ($\{\omega_i\}_\mathrm{shuffled} \to \mathrm{shuffle} \left( \{\omega_i\}_\mathrm{original} \right)$ or by changing the frequency of only a single unit to a new value ($\omega_{i, \mathrm{original}} \to \omega_\mathrm{\mathrm{new}}$). We show a transition to phase synchronization with increasing coupling strength and switching from short-range to long-range connections, as expected, and find large sample-to-sample (STS) fluctuations (Fig. \ref{fig:transition_sync}) 
We study two coupling schemes, Watts-Strogatz (WS, small-world) and distance-dependent (DD).
We consider ensembles as collections of networks with fixed coupling strength $\epsilon$ and topology (fixed rewiring probability $p$ or locality parameter $\alpha$) but distinct realizations of the natural frequencies $\{\omega_i\}$ \cite{carlson2011sample}.  Each ensemble in the figure contains $501$ samples (realizations). We present the results using the mean degree of phase synchronization $R$ for each realization, and the gap $\Delta \walrus R_\mathrm{max} - R_\mathrm{min}$ between the most and least phase synchronized realizations in each ensemble. \kl{The gap $\Delta$ is chosen simply to illustrate the wide range of $R$ values clearly, and we remark that very similar curves are observed by using the standard deviation over samples.}
In Fig. \ref{fig:transition_sync}, thicker lines represent an ``original" sequence of frequencies $\{\omega\}_\mathrm{original}$, from which other realizations (light lines) are created by shuffling all frequencies or by changing one unit's frequency to a new value $\omega_\mathrm{new} = 3.0$. Each sample is a different dynamical system, and has a different transition to phase synchronization, which occurs at different values of $\epsilon$, $p$, or $\alpha$, and with a different profile (some can have a region of a small region of desynchronization while others do not, for instance). 

This means that changing samples can lead to large changes in the behavior of the system. Let us firstly study the transitions induced by increasing the coupling strength $\epsilon$ for four representative types of networks (panels (a)-(d)), characterized by four specific values of rewiring probability $p$ and locality parameter $\alpha$.
In the red curves, networks are dominated by long-range connections, with $p = 1$ (random) and $\alpha = 0$ (all-to-all)  and have a complete transition to phase synchronization (reaching $R \sim 1$), with sample-to-sample fluctuations (measured by $\Delta$) increasing during the transition and returning to zero after.
The all-to-all (mean-field) case is the finite-size version of the system originally studied by Kuramoto  \cite{kuramoto1984chemical}, and the critical $\epsilon$ values , when the transition occurs in each sample, are close to the $\epsilon_c = \frac{2}{g(0)\pi} = \frac{2\sqrt{2}}{\sqrt{\pi}} \approx 1.596$ predicted in the thermodynamic (infinite network size) limit. Its finite-size scaling properties have also been studied in \cite{hong2007entrainment}. It is worth mentioning that this parallel between random networks and mean-field networks, which have similar phenomenology, has been described in other works, with random topologies having the same scaling exponents, and similar transition, than the mean-field (all-to-all) topology \cite{skardal2020higher, hong2013link}. 

In the green curves ($p = 0.19684$ and $\alpha = 1.538463$) some connections have been rewired in the Watts-Strogatz networks, and weights redistributed for distance-dependent networks, from long-range to effectively short-range connections. On average, phase synchronization $R$ decreases, though still remaining high. Some samples of WS networks also start to display regions of desynchronization: after the initial transition to high $R$, a further increase in $\epsilon$ can desynchronize them (visible in panels (a) and (c), for $\epsilon$ roughly in $[6,7]$). Therefore, the huge changes in $R$ ($\Delta \sim 0.99$) due to changing samples can be attributed to two effects: the difference in their critical coupling strength (when the transition begins) \cite{hong2006anomalous}, and also in their different post-transition behaviors (such as the desynchronization gaps that occur at differing intervals of $\epsilon$.)

In the purple curves ($p = 0.08733$ and $\alpha=1.76923$), even more short-range connections become present. Phase synchronization $R$ on average decreases, and the fluctuations $\Delta$ remain high, though occur more evenly spread across samples.
Finally, for cyan curves ($p = 0$, two-nearest-neighbor chains and $\alpha = 3$, close to nearest-neighbor chains) all connections are short-range. Their phase synchronization is much smaller, and they do not reach a high degree of phase synchronization for any value $\epsilon$ we tested.
These networks with short-range connections still have STS fluctuations in $\meanR$, though at a smaller degree than the others. 

Returning to frequency synchronization, we mention that for weak coupling strengths (roughly below $\epsilon \approx 3$), most of the samples in any ensemble are not frequency synchronized (see Fig. S3). Above this value, frequency synchronization becomes more common, especially for networks with more long-range connections, such that for sufficiently high coupling all samples become frequency synchronized. This is not the case for networks with mostly short-range connections ($p \lessapprox 0.01$), in which some samples do not reach frequency synchronization even despite strong coupling. The presence of frequency synchronization in the short-range networks is consistent with the literature \cite{strogatz1988phase, acebron2005kuramoto} showing that frequency synchronization in first-nearest-neighbor chains is possible for sufficiently high $\epsilon$ in strictly finite systems. There are therefore also sample-to-sample fluctuations in the frequency synchronization of Kuramoto networks. They occur similarly to the fluctuations in phase synchronization, but are somewhat harder to visualize and have a less interesting dependence on parameters, justifying our focus on phase synchronization in this paper.

We now move to the topology induced transitions, which occur by switching from short-range to long-range connections (varying $p$ and $\alpha$) while keeping the coupling strength $\epsilon$ fixed (Figs. (e-h)). A similar scenario with a continuous transition to phase synchronization occurs, induced by changing either $p$ or $\alpha$. The STS fluctuations $\Delta$ increase during the transition, reaching significant values for both shuffled realizations and single-unit changes. The nearest-neighbor networks show some STS fluctuations, while the long-range dominated ones (random or mean-field) do not. We note here that the transition for WS occurs at $p \sim 0.1$, so we plot the figures in log scale to show the full transition to synchronization. This transition was already reported for WS networks in \cite{hong2002synchronization}, but the authors used a linear scale for $p$ and missed the full details of the transition that we see here, especially the sample-to-sample fluctuations; for distance-dependent powerlaw networks, a transition in phase and frequency was reported in \cite{rogers1996phasetransitions}. None of these references studied the sample-to-sample fluctuations, however. 

We conclude that either shuffling or changing a single unit can significantly alter the behavior of these systems, leading to large sample-to-sample fluctuations, in some cases over a very large range of parameters. This is particularly true for WS networks, reaching $\Delta \sim 0.99$, close to the maximum possible value of $\Delta = 1.0$. The distance-dependent networks have weaker fluctuations, though still significant, reaching up to $\Delta \sim 0.7$.

\begin{figure*}[htb]
    \centering
    \includegraphics[width=0.9\textwidth]{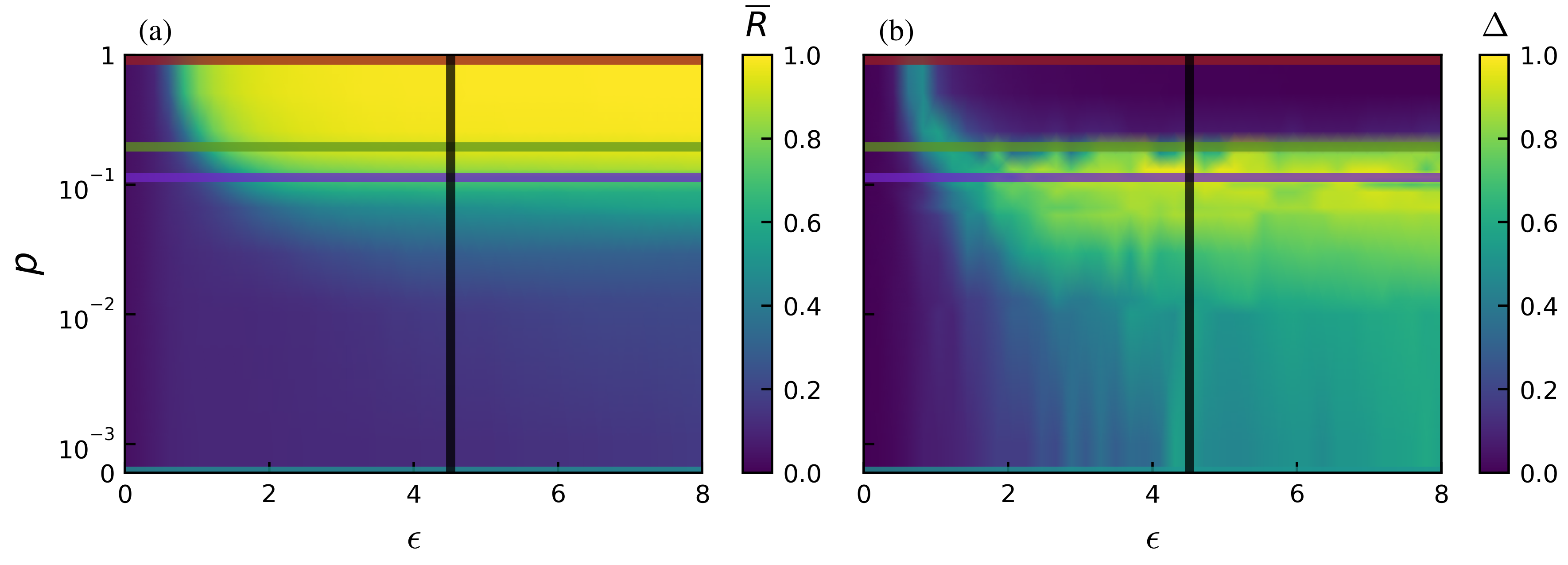}
    \caption{\textbf{Sample-to-sample fluctuations increase around the transition regions to phase synchronization.} The surface on the left shows the average degree of phase synchronization $\overline{R}$ across the ensemble ($1000$ realizations of shuffled natural frequencies). The region of high phase synchronization is clearly seen for sufficiently high coupling strength $\epsilon$ and rewiring probability $p$. The colored lines correspond to the parameter values shown in Fig. \ref{fig:transition_sync}. The right panel displays $\Delta$, the difference between the most and least synchronized realizations for each pair $(p, \epsilon$), and we see that the fluctuations from sample to sample increase during the transitions to phase synchronization. The green and purple curves remain close to the region of transition for all $\epsilon \gtrapprox 1$, such that their fluctuations do not decrease with an increase in $\epsilon$. The figure uses Gouraud interpolation to ease visualization by smoothing the curves with a linear interpolation.}
    \label{fig:Rsurface}
\end{figure*}

The networks with intermediate $p$ or $\alpha$ and the short-range networks have STS fluctuations even for high $\epsilon$. This is consistent with the known increase in the fluctuations near a phase transition because the networks with these parameters remain close to the topology-induced transition.
This is illustrated for WS networks in Fig. \ref{fig:Rsurface}, which shows, in the $p \text{---} \epsilon$ parameter space, the average $\overline{\meanR}$ phase synchronization $R$ across samples on the first panel and the sample-to-sample fluctuations measured by $\Delta$ on the second panel. Figure \ref{fig:Rsurface} provides a comprehensive view on both the coupling strength and the topology-induced transitions. The samples are realized here as shuffles, though a similar figure would be obtained by changing one unit. There is a single region of phase synchronization for sufficiently high coupling strength $\epsilon$ and rewiring probability $p$ (panel (a)). Around the borders of this region, where the system is transitioning, the sample-to-sample fluctuations are much higher (panel (b)). It then becomes clear that the intermediate networks (green and purple lines), are near the topology-induced transition (for instance, black line) for all $\epsilon \gtrapprox 1$. As $\epsilon$ is increased, they remain near this $p$-transition, and so their sample-to-sample fluctuations do not decrease. 
For the regular networks, we first note that the $p$-axis is shown in a log scale, such that these networks, with $p = 0$, are still relatively close to the transition at $p_c \approx 0.1$, and thus they also present significant STS fluctuations.

Figure \ref{fig:Rsurface} also illustrates the existence of two qualitatively different types of transitions: one induced by increasing coupling strength (for sufficiently high $p$), another induced by increasing $p$ (for sufficiently high $\epsilon$). The difference between both is in their starting points. Both are globally phase desynchronized, but in the former (red, green and purple lines), the weak coupling strength regimes have mostly uncorrelated oscillators, with no discernible structures in the phases or even synchronization in the frequencies. In the latter (black line), there are shorter-range structures with frequency synchronization for most samples.

\subsection{Unpredictability of dynamical malleability}
For Watts-Strogatz networks, samples can be generated by resampling the topology instead of changing the natural frequencies. Since they are generated by a random rewiring process, different realizations generate different networks (there is link-disorder \cite{hong2013link}). Different samples can therefore also be generated by resampling the network while keeping the natural frequencies fixed. This generates a profile of sample-to-sample fluctuations similar to that shown in Figure \ref{fig:transition_sync}(e), where the network was fixed and the natural frequencies were changed (see Fig. S1 for details). 

\begin{figure*}[htb]
    \centering
    \includegraphics[width=0.9\textwidth]{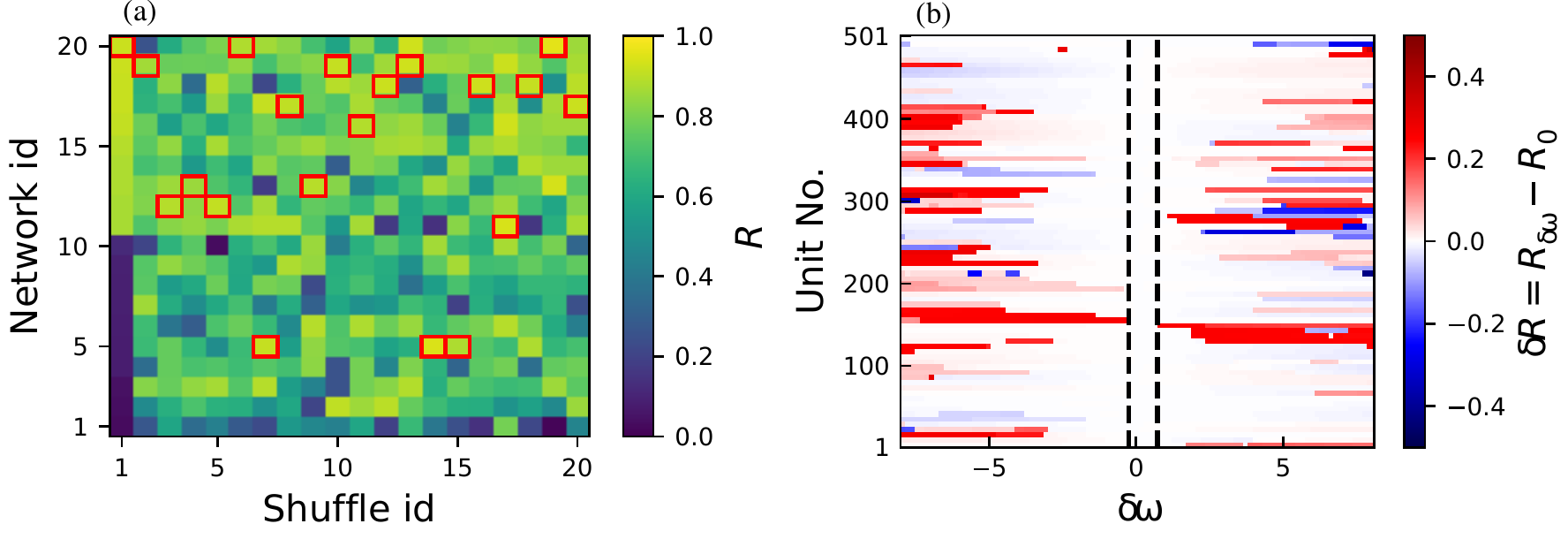}
    \caption{\textbf{Fluctuations in dynamically malleable systems are unpredictable}. Panel (a) shows the average phase synchronization $\meanR$ for fixed coupling strength $\epsilon = 4.51282$ and rewiring probability $p = 0.08733$ for different $20$ shuffles of the natural frequencies $\{\omega_i\}$ and samples of networks generated by the Watts-Strogatz algorithm. For ease of visualization, the networks are ordered such that the highest Network ids correspond to higher synchronization for $\mathrm{Shuffle\;id} = 1$. For each shuffle, the network with the highest $\meanR$ is marked with a red rectangle. We thus see that no network synchronizes more for all shuffles: $\meanR$ is a function of both the specific frequency and topology samples. Panel (b) shows the changes $\delta R$ in the phase synchronization $R$ when the natural frequency of each unit is changed by an amount $\delta \omega$, such that $\omega_i \to \omega_i + \delta\omega$. Other parameters are fixed, in particular $p = 0.1145$ and $\epsilon = 4.51282$. There is a rough threshold (indicated by the black dashed lines), below which changing $\omega_i$ does not significantly alter $R$ ($\delta R < 0.1$ for the figure). Furthermore, changing the frequency does not have a monotonic impact on the change in $R$: small alterations in $\omega_i$, above the threshold, can have the same impact on $R$ as bigger alterations. 
    }
    \label{fig:complexxsensitivity}
\end{figure*}

Now, we wish to illustrate that no network, or natural frequency sequence, are alone responsible for leading to more, of less, synchronized states. Instead, the samples depend sensitively on both, especially in the region of large STS fluctuations. Figure \ref{fig:complexxsensitivity}(a) shows the degree of phase synchronization for different realizations of the networks and different shuffles of the natural frequencies, all for $\epsilon = 4.51282$ and $p = 0.08733$ with fixed initial conditions. To aid the visualization, red rectangles indicate the network with the largest $R$ for each shuffle. No network synchronizes more (or less) for any sequence of natural frequencies; and no sequence of natural frequencies synchronizes more for any network. Furthermore, if the $\epsilon$, $p$, or initial condition is changed, the whole profile of the figure also changes.
Another way to illustrate the complex sensitivity in the region of high sample-to-sample fluctuations is by now fixing the network, and changing the frequency of a single unit by an amount $\delta \omega$. Fig. \ref{fig:complexxsensitivity}(b)) illustrates the change $\delta \meanR$ in the phase synchronization, compared to the synchronization of the "original" ($\delta \omega = 0$) frequency realization. There is a rough threshold, at $| \delta \omega | \gtrapprox 0.1$, below which the perturbations in one unit do not significantly affect the network's phase synchronization. Above this threshold, however, large changes occur. They are asymmetric on $\delta\omega$ and occur non-monotonically (increasing $|\delta\omega|$ does not necessarily to bigger changes). This complicated pattern we observe could make the design and control of these systems quite difficult in practice.

\subsection{Ratio of short to long-range connections}
The rewiring of connections, in WS networks, or the redistribution of weights in DD networks, from short-range to long-range leads to a transition towards globally phase-synchronized regimes. During this transition, the sample-to-sample fluctuations peak for some ratio of short-range to long-range connections. To quantify this ratio, we first define the short-range connections to/from a node $i$ as all existing connections to/from other nodes $j$ within an edge distance $d$ (with index $j \in [i-d, i+d]$), with $d$ being the range of short connections ($d = 2$ here). For WS networks, we calculate the average degree (number of connections) for short-range ($K_s)$ and long-range connections ($K_l$). For DD networks, we define an analogous measure of topological influence. They are
\begin{align}
    K_{s} & \walrus \frac{2}{\eta(\alpha)} \sum\limits_{j=1}^{d} \frac{1}{j^\alpha} \\
    K_{l} & \walrus \frac{2}{\eta(\alpha)} \sum\limits_{j=d+1}^{N^\prime} \frac{1}{j^\alpha} .
\end{align}

Note that due to the symmetry of the DD networks, nodes share the same value of $K_s$ and of $K_l$. The ratio $\kappa$ of short-range to long-range connections is then defined as:
\begin{equation}
    \kappa \coloneqq \frac{ K_\mathrm{s} - K_\mathrm{l} }{ K_\mathrm{s} + K_\mathrm{l} },
\end{equation}
so that $\kappa = 1$ if only short-range connections exist, and $\kappa = -1$ if only long-range connections exist, with intermediate cases in between. In WS networks, the number of connections is $K = kN$ ($k$ being the amount of neighbors of each node), with the number of long connections approximately $K_l = pK$ and short-range approximately $K_s = (1-p)K$. Therefore, the ratio $\kappa$ can be easily calculated to be approximately $\kappa = 1 - 2p$. For DD networks, the ratio $\kappa$ is given as 
\begin{equation}
    \kappa = \frac{ \sum_{i=1}^d i^{-\alpha} - \sum_{i=d+1}^{N^\prime} i^{-\alpha} }{\sum_{i=1}^{N^\prime} i^{-\alpha}  }.
\end{equation}

\begin{figure*}[htb]
    \centering
    \includegraphics[width=0.75\textwidth]{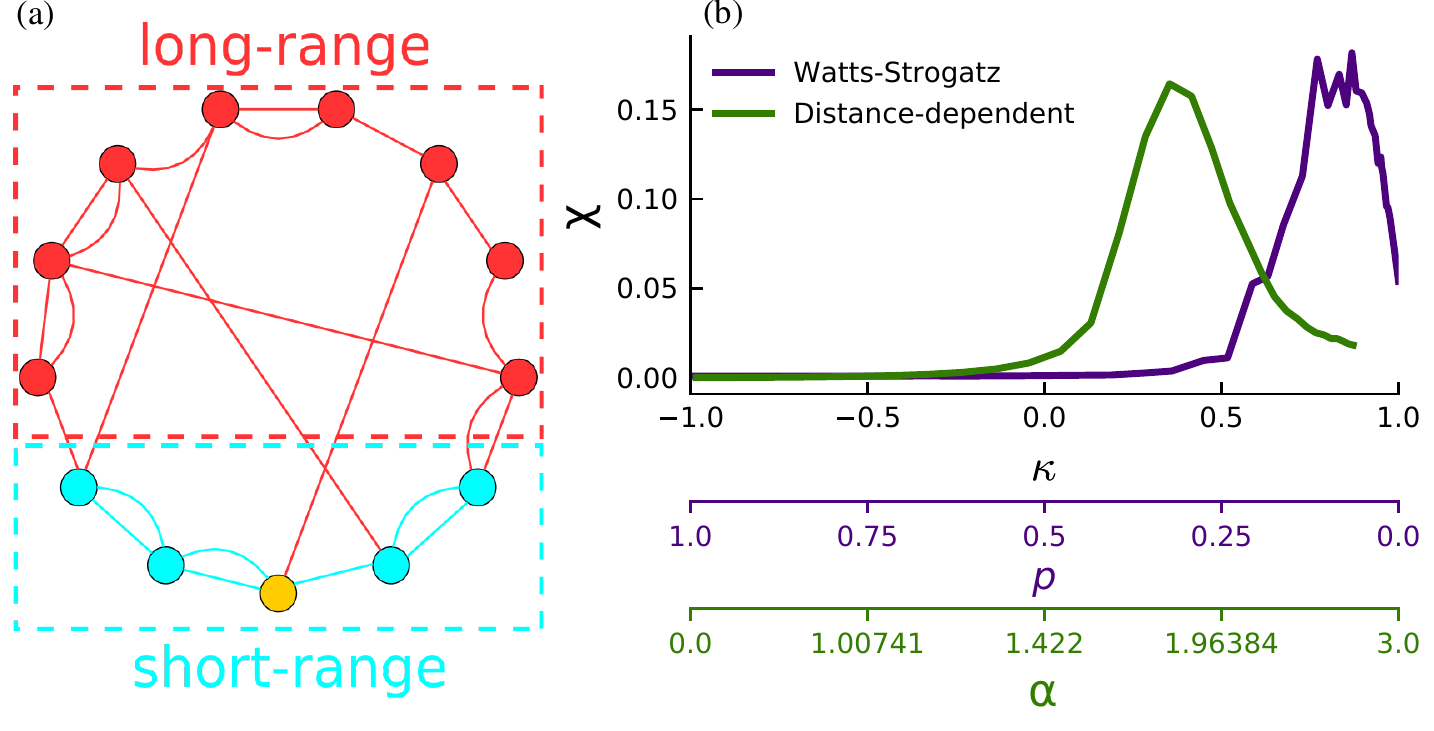}
    \caption{\textbf{Sample-to-sample fluctuations peak within a narrow interval in the relation of short-range to long-range connections}. Panel (a) illustrates the short-range (blue) and long-range (red) connections from the yellow unit for $d = 2$. Panel (b) shows the sample-to-sample fluctuations in the phase synchronization measured as the standard deviation $\chi$ of the distribution function of $\meanR$ against the ratio $\kappa$ of short-range to long-range connections calculated for several distinct topologies $p$ and $\alpha$. The green curve corresponds to the distance-dependent networks, with $\epsilon = 6.46154$ and 501 realizations per $\alpha$; purple corresponds to Watts-Strogatz networks, $\epsilon = 4.51282$ and 1501 realizations per $p$. The bottom axis show the values of $p$ and $\alpha$ for the respective ticks in $\kappa$ (note the spacing for $\alpha$ is not linear). }
    \label{fig:ratioshortlongrange}
\end{figure*}

Figure \ref{fig:ratioshortlongrange} shows this ratio $\kappa$ calculated for the same setup of Fig. \ref{fig:transition_sync}(e) and (f), shuffling natural frequencies with fixed coupling strength and changing $p$ or $\alpha$. The STS fluctuations are measured here by standard deviation $\chi$ across the samples, instead of $\Delta$. The former makes the figure clearer, but the same analysis also works using $\Delta$. A remark when comparing with Fig. \ref{fig:transition_sync} is that the two measures may peak at slightly different values of $p$ or $\alpha$.  For both types of networks, the STS fluctuations peak when there is a relatively small number of long-range connections present in a short-range-dominated network. It is more extreme for WS, as the ratios are closer to $1$ than in the DD networks. This discrepancy in the ratios leading to higher STS fluctuations shows that $\kappa$ is not an universal feature for any topology, but can still be important to understand the behavior.

\subsection{Multistability}
So far, we have changed natural frequencies while keeping initial conditions fixed. Now we invert this, and shuffle initial conditions to study the system's multistability. We continue examining phase synchronization (PS) $\meanR$, though we know that $\meanR$ is only a rough measure of multistability. Being a mean value, the same $\meanR$ could represent different attractors. Therefore, the number of attractors estimated based on $\meanR$ can only be considered as a lower bound. To remedy this, we also verified the findings by comparing several other features of the dynamics. amhese included the standard deviation of PS in time, the PS between each unit and its neighbors, the PS between sections of $100$ units, the time-averaged instantaneous frequencies $\dot{\theta}_i$ of units, and the standard deviation, inter-quartile interval and gap between the unit's instantaneous frequencies. Realizations with unique values of all these features were considered as a distinct attractor. The number of such attractors agrees qualitatively with the dispersion we see in $\meanR$, increasing during the transition.

\begin{figure}[htb]
    \centering
    \includegraphics[width=1.0\columnwidth]{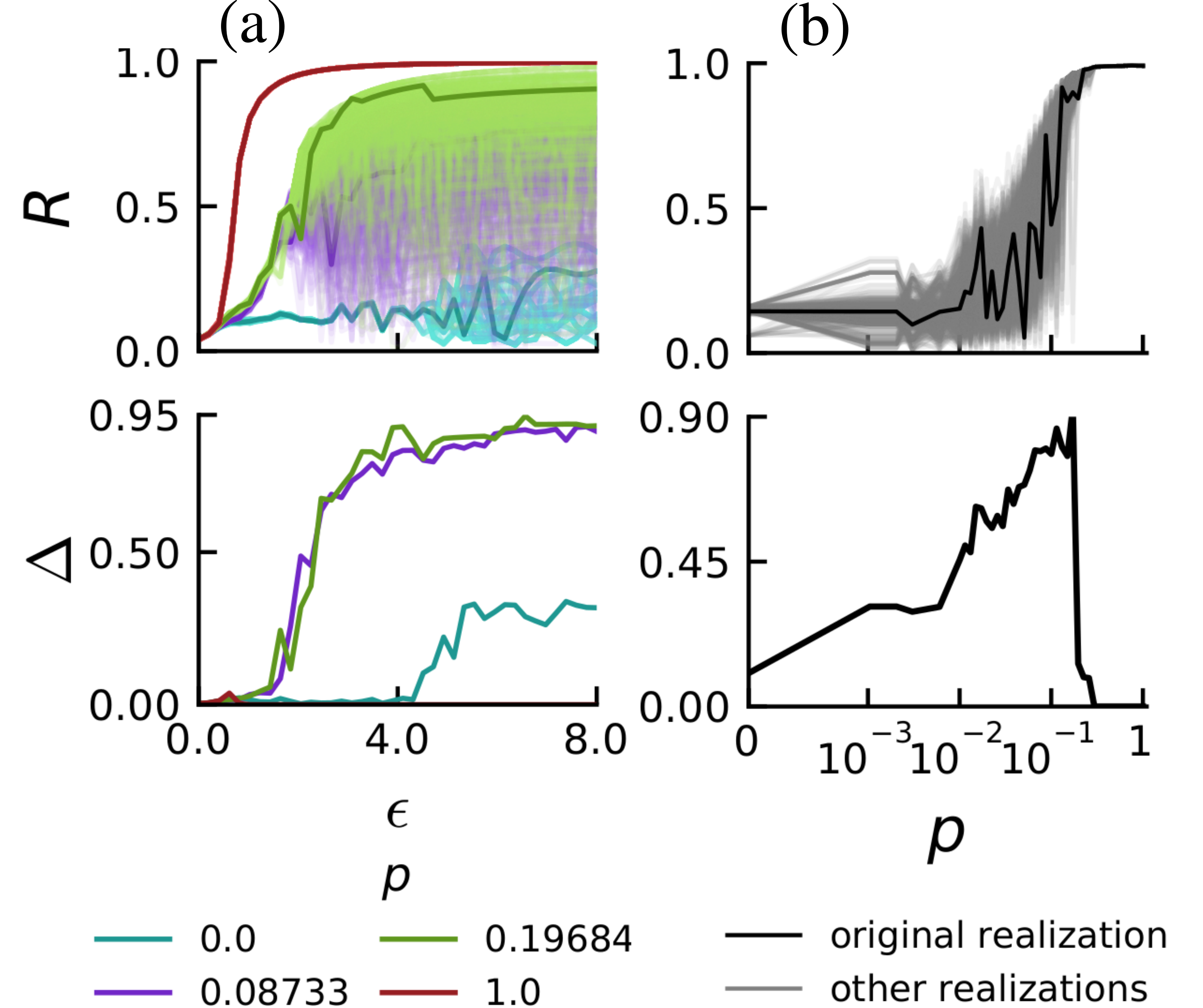}
    \caption{\textbf{Multistability in WS networks}. Phase synchronization and its dispersion for $501$ different shuffles (thinner lines) of the initial conditions, taken from the original initial conditions (thicker lines) used throughout the rest of this work. All other parameters are fixed, including the natural frequencies as the original frequency distribution. The coupling strength $\epsilon$ (left panel) and rewiring probabilities $p$ (right panel) are the same ones used for WS networks in Fig. \ref{fig:transition_sync}. The multistable behavior is thus very similar to what we observed before by changing the frequencies (Fig. \ref{fig:transition_sync}(a) and (e)), and so shuffling the initial conditions for this network also leads to large fluctuations in the phase synchronization. }
    \label{fig:multistability}
\end{figure}

The phase synchronization is thus shown in Fig. \ref{fig:multistability}. Random networks ($p = 1$, red) are multistable only during their transition to phase synchronization. Intermediate networks ($p = 0.19684$, green; $p = 0.08733$, purple) have a high degree of multistability, meaning coexistence of several attractors, with very distinct degrees of phase synchronization. No shuffle of the initial conditions leads here to the same attractor, so the system has at least $501$ attractors, the number of different realizations tested.
The 2-nearest-neighbor lattice has significant multistability for $\epsilon \gtrapprox 4$. This is consistent with the literature for 1-nearest-neighbors, in which multistability occurs after the transition to phase-locking \cite{tilles2011multistable}. 

This multistability can enhance the sensitivity of the system to parameter changes, and help to explain the large fluctuations we observe. In this case, a parameter change needs only to change the boundaries of the basins of attraction for the same initial condition to land on a completely different attractor. Attractors do not have to be necessarily drastically changed for the large STS fluctuations to be observed. However, multistability is not in principle required for STS fluctuations; in fact, the distance-dependent networks appear to be monostable (not shown), though they are malleable. 

\subsection{Distributions of samples}
\begin{figure*}[htb]
    \centering
    \includegraphics[width=0.95\linewidth]{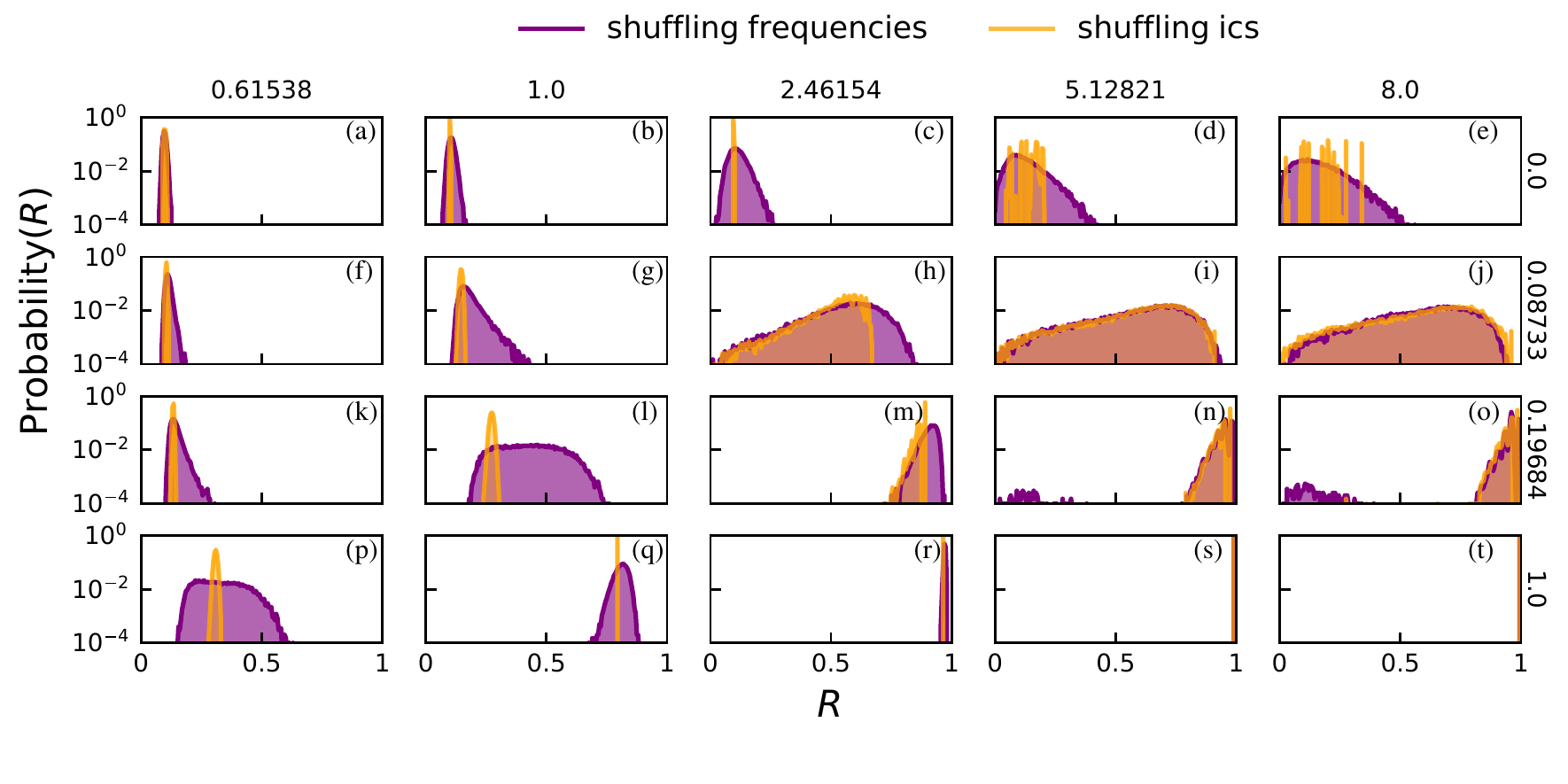}
    \caption{\textbf{Distributions of $\meanR$ due to shuffling frequencies or initial conditions.} Each panel contains the distribution of the mean degree of phase synchronization $\meanR$ across 20000 shuffles of natural frequencies (in purple) or initial conditions (orange) for Watts-Strogatz networks. The rewiring probabilities $p$ are indicated on the right of each row, and are the same as used in Fig. \ref{fig:transition_sync}(a); the coupling strengths $\epsilon$ are indicated on the top of each column. Bin size is $0.005$, and the probability for each bin is calculated as the occupation of the bin divided by the total occupation across all bins, and is shown in logarithmic scale.}
    \label{fig:Rdistributions}
\end{figure*}

As we have seen, shuffling initial conditions can also generate realizations with widely different dynamics, similarly to shuffling natural frequencies. But the two methods to create an ensemble of samples have different effects, and can generate samples with distinct distributions. As shown in Fig. \ref{fig:Rdistributions} for Watts-Strogatz networks, shuffling frequencies leads usually to a broader, and smoother, distribution of $\meanR$. This increased broadness shows that new attractors are indeed created by shuffling the frequencies, so that multistability itself cannot account for the sample-to-sample fluctuations we discussed previously. Furthermore, the transitions to phase synchronization occur through an increase in the distribution's average. The accompanying increase in the width of the distribution shows the increase of sample-to-sample fluctuations, which go to zero only for long-range networks ($p = 1$). 

Specifically, the distributions for the two-nearest-neighbor lattice ($p = 0$, panels (a)-(e)) are quite different: shuffling frequencies leads to a smooth distribution, whose average shifts to the right as $\epsilon$ is increased; for shuffling initial conditions there is also a slight increase in the distribution's average as $\epsilon$ is increased, but the distribution itself is dominated by several peaks. 
For intermediate networks ($p = 0.08733$ and $p = 0.19684$, panels (f)-(o)), the skewness of the distribution becomes negative, and shuffling initial conditions has a smoother behavior, more similar to shuffling frequencies. Interestingly, the distribution can be bi-modal, with the two modes being separated on either extreme of $\meanR$ (panels (n) and (o)). 
For $p = 1$ (random network), the two first coupling strengths (panels (p)-(q)) occur during the narrow interval of significant STS fluctuations, during the transition to phase synchronization. Soon after $\epsilon > \epsilon_c \approx 1.6$, the distribution becomes extremely narrow. 

It is worth mentioning that very similar distributions are obtained if, instead of shuffling the frequencies or initial conditions, we re-sample them from the distribution (i.e. change the seed in the random number generator).

\section{Discussions and conclusions}\label{sec:discussion_conclusion}
In this work we have studied ring networks of Kuramoto phase oscillators near the transition to phase synchronization. Our main findings contribute to a deeper understanding of the high sensitivity of a network's dynamics with respect to changes in the natural frequency of its oscillators, leading to a large variety of different paths to phase synchronization, i.e. large dynamical malleability of the network. The choice of networks of Kuramoto oscillators for the analysis is based on the fact that this model is a generic representation of weakly coupled oscillators in which the amplitude dynamics can be neglected. The networks considered are heterogeneous in their units - all nodes have different parameters, here their natural frequencies - and have different topologies - Watts-Strogatz (small world) topology or a distance-dependent topology. 
Our most striking observation is that in such networks even a small perturbation to the parameter of a single unit can drastically change the dynamics of the whole network. Considering each change to the frequency of a single unit as a particular sample for the network's dynamics, we find large sample-to-sample (STS) fluctuations in the path to phase synchronization as has been reported previously for globally and randomly coupled networks of Kuramoto oscillators \cite{hong2007entrainment, peter2018transition}. These large fluctuations are manifested in the dynamics of the phase patterns in the network, which are different for each choice of a single new frequency, as well as for shuffling all the frequencies.  

In general, increasing the coupling strength leads to a transition from uncorrelated phases to the existence of patterns in the phases, as it makes the units' instantaneous frequencies synchronize (cf. Fig. S3). Which phase patterns are present depends not only on the distribution of the natural frequencies but also on the network topology, i.e. its ratio of short to long-range connections. Replacing shorter-range connections with longer-range ones leads to networks with longer-range patterns, that are more globally phase synchronized. Therefore we observe in fact two qualitatively different transitions to phase synchronization, with distinct initial, desynchronized regimes: in the coupling strength transition, for weak coupling strengths the oscillators are very weakly correlated, such that their phases do not form patterns and their instantaneous frequencies are not synchronized; in the topology transition, too few long-range connections make the system also phase desynchronized, but it nevertheless has short-range patterns and is generally frequency synchronized, as long as the coupling strength is sufficiently large.

Because of the coexistence of the two transitions, the dynamical malleability remains high over a wide range in parameter space which is spanned by coupling strength $\epsilon$ and rewiring probability $p$ for Watts-Strogatz networks. Networks with an intermediate amount of long-range connections are highly malleable for any coupling strength $\epsilon$ we tested, starting at a relatively weak value of coupling strength ($\epsilon \sim 1$). This is because they always stay close to the topology-induced-transition, even though they are far from the coupling-strength-induced-transition. 
For these intermediate networks, an increase in coupling strength can also significantly alter the phase synchronization even after the networks have transitioned to phase synchronization. They can go from highly phase synchronized, to desynchronized, and then return to synchronized (as seen in Fig. \ref{fig:transition_sync}(a),(c)). This is a known behavior in systems with finite size, and can be observed in phase transitions in simpler systems, such as cellular automata. It means that the systems can be highly sensitive not only to changes in the units' parameters, but also in the global coupling strength, even after the transition has taken place.

It is worth noting that Watts-Strogatz networks are typically more malleable than the distance-dependent ones. One reason for this is that WS networks have an additional parameter disorder due to the random switching of connections (link-disorder \cite{hong2013link}): because these systems are finite, different realizations, with the same value of $p$, can lead to different networks. This is not the case for DD networks, for which each $\alpha$ defines a unique network, and thus have no topological disorder. Indeed, fluctuations of different topology realizations, with fixed initial conditions and natural frequencies, have a similar phenomenology to fluctuations due to frequency realizations (Fig. S1), and are thus also a source of dynamical malleability for the system.  

The transitions to phase synchronization correspond to non-equilibrium phase transitions \cite{kuramoto1984chemical, peter2018transition}, such that we can connect the dynamical malleability we analyze with the sample-to-sample (STS) fluctuations well-known in statistical physics. There, it is usually described in finite systems that different samples have different statistical properties, which lead to distinct phase transitions - transitions are usually said to be shifted between samples \cite{sornette2006critical, hong2007finitesizescalingpre, hong2007entrainment}. In this work, we show that there are many details beyond this, and also highlight that the STS fluctuations can (i) have very large magnitudes and (ii) depend in a complicated manner on the parameters of the system. These points are still under-explored in the literature, but can have a significant impact in systems with transitions to synchronization.
We have shown the large magnitude of the STS fluctuations in a variety of manners: (i) changing the natural frequencies of single units separately or shuffling all in the network; (ii) changing the realization of the network topology and (iii) for two different types of topology (Watts-Strogatz and distance-dependent). In all cases, the global state of the network can be drastically disturbed, moving from desynchronized to synchronized, for instance. We have also explored the complicated dependence of the behavior in different ways. We have pointed out that the distinctions between samples are not just in the start of their transitions, like the usual case considered in statistical physics. After the initial transition to phase synchronization, samples can have large dips of desynchronization as the coupling strength increases. These dips result in large fluctuations between samples since, for a same coupling strength, one sample may dip to a desynchronized state while another sample remains highly synchronized (see the thicker green curve in Fig. \ref{fig:transition_sync}(a) for an example).

Furthermore, we have shown that the behavior of the systems is a complicated function involving the coupling strength, topology, and natural frequencies together: we have not found a sequence of frequencies, or a specific network realization, that always leads to more (or less) synchronized networks (Fig. \ref{fig:complexxsensitivity}(a)). Looking at the pair of sequence of frequencies and network realization in the most phase synchronized sample, for instance, we also notice (not shown) that they change depending on initial condition and coupling strength. Changes in the natural frequency of single units also lead to non-monotonic changes in the synchronization. In other words,  varying the change in frequency can increase or decrease the synchronization level, depending on the chosen unit, coupling strength and topology (Fig. \ref{fig:complexxsensitivity}(b)). As a consequence, we are unable to identify a specific unit, or magnitude of perturbation, that is always responsible for the greatest disruption. It is important to mention, though, that there seems to be a minimum magnitude of perturbation needed to disrupt the system significantly: the change $\delta \omega$ in one unit's natural frequency needs to be bigger than around $\delta \omega \gtrapprox 0.1$. 

Additionally, we find that changing the initial conditions leads to a similar phenomenology compared to changing the natural frequencies for WS networks. That is, they are \textit{multistable} in a wide parameter range, being monostable apparently only for networks close to or at the random regime. This, we believe, is a novel result in the literature. As a consequence, the coexistence of multiple attractors for a given set of the unit's parameters and topology can enhance the STS fluctuations. Moreover, changes to the natural frequencies, which lead to  changes in the attractors and/or in their basins of attraction, will result in a different dynamics though the initial conditions have been fixed. However, multistability is not required for dynamical malleability due to changing frequencies, as we see in all-to-all and distance-dependent networks, which do not appear to be multistable but still have large dynamical malleability. Multistability could have an even stronger impact in the malleability if the basins of attraction were complexly interwoven. Then, even very small changes could lead to significant STS fluctuations. But this does not appear to be the case in any of the networks we studied, which seem to have smooth basin boundaries (Fig. S4) - it is thus noteworthy that the already high dynamical malleability we have described can occur even with smooth boundaries.

Changing initial conditions, however, leads to different statistics than changing frequencies. We have observed that the distribution of the degree of phase synchronization for several samples of natural frequencies and initial conditions, generated by either shuffling or re-sampling them, is smoother for frequencies' changes, indicating a higher number of possible attractors. The distribution is also broader, meaning fluctuations are stronger due to changing frequencies. This shows that changing frequencies creates new attractors, and the large STS fluctuations are not simply a consequence of the system's multistability. Interestingly also, the shape and skewness of the distributions change for different topologies as $p$ is increased. The distributions are not Gaussian, which is inconsistent with the assumptions made in other works \cite{hong2007entrainment, hong2013link}. In those works, the authors argue that the fluctuations must be normally distributed for sufficiently large networks and many samples due to the central limit theorem. This inconsistency is likely generated by the finite size of the networks studied here. Even in all-to-all networks, in which there is no topological disorder, the distributions are not Gaussian for $N=501$, but start approaching Gaussian distributions as $N$ is increased to $5000$. 

The size of the networks influences, thus, the distributions. It also influences the magnitude of dynamical malleability and the interval of parameters in which it occurs. The results presented in the paper are for networks with $N = 501$ oscillators. Scaling analysis (Fig. S2) reveals that the intervals of high malleability decrease with the size $N$, as expected from the theory - for instance, authors in \cite{hong2007entrainment} describe the range of $\epsilon$ for high malleability as scaling as $N^{-2/5}$ for all-to-all networks. 
For the WS networks, malleability is still significant for even up $N = 5000$ oscillators. Moreover, the maximum magnitude of the fluctuations does not decrease significantly, and networks with $N = 5000$ can still reach $\Delta = 0.9$. 
This suggests that the malleability gets restricted to a smaller region in parameter space, but might not decrease significantly in magnitude there for bigger networks. In the limit of infinite-size networks, it would get restricted to a single line, defining the two transitions to phase synchronization, and remain non-zero there. This is consistent with a study in global networks of Kuramoto oscillators, where this behavior was observed \cite{hong2006anomalous}. In fact, systems of this type are well-known in phase transitions with quenched disorder (heterogeneous parameters), where they are said to be non-self-averaging \cite{wiseman1995lack}. For these systems, the STS fluctuations (i.e. dynamical malleability) gets restricted to a smaller region of parameter space, but nevertheless remains finite there \cite{wiseman1995lack}. This region is the point, for one parameter spaces, or manifold, for higher-dimensional parameter spaces, at which the phase transition occurs - they are the critical point, or critical manifold. 
In any case, networks of $N = 5000$ units can be regarded as rather large in several real-world applications \cite{peter2018transition}, so the STS fluctuations we describe here occur for a significant range of system sizes including those of several applications. 

We remark again that the phenomenology we describe also occurs for wide ranges of topology and coupling strength values, for distinct frequency distributions, such as Cauchy-Lorentz instead of Gaussian (not shown), and for other dynamical models. For instance, previous work on spiking neuron networks has revealed a very similar phenomenology \cite{budzinski2020synchronization}. Unpublished results also show similar behavior for cellular automata of spiking units. We have observed (not shown) similar behavior in small-world networks generated by adding long-range connections and keeping the short-range ones \cite{newman1999scaling}. Other works have also observed dynamical malleability in Kuramoto oscillators coupled under both human-connectome structural networks and hierarchical-modular networks \cite{buendia2022the, villegas2014frustrated}. Additionally, of course, the theory of phase transitions and, consequently, of sample-to-sample fluctuations, is known to apply for a variety of distinct systems. Furthermore, malleability does not seem to be limited to a specific parameter type: in the Kuramoto network described here, changes made to any parameters, both natural frequencies and topology, can lead to large fluctuations.

This leads to an interesting question: is dynamical malleability good or bad for systems? 
On one hand, large fluctuations can be undesired. For instance, a large fluctuation could take power grids from a phase synchronized regime to a desynchronized one, and lead to blackouts.
On the other hand, fluctuations can be desired due to the increased flexibility of the systems. They could be a useful mechanism for adaptation, learning or memory formation in neural circuits. More specifically, an important property of the brain is that it can separately process information from different types of input in segregated areas, and then integrate them all into a unified representation \cite{tononi1994a, deco2015rethinking}, a process that is thought to be important in consciousness \cite{tononi1998consciousness}. For this reason, Tononi and colleagues conjectured that the brain needs to have an optimal balance between segregation and integration of areas \cite{tononi1994a}. In this optimal balance, synchronization between different regions in the brain 
needs to fluctuate considerably, ranging from low synchronization to high synchronization \cite{fingelkurts2006timing}. Therefore, having a large dynamical malleability can be an advantageous feature, allowing for this high variability to be achieved through small changes in the neurons, e.g. their firing rate, or their connections. There is also an interesting evidence for this in \cite{li2009burst}, which reported that high-frequency firing of neurons can drive changes in the global brain state. Future research on dynamical malleability can thus be dedicated to understand ways to quench or to explore the fluctuations, using the framework we establish here.

An interesting line of research is also opened here by considering the effect of noise or time-dependent forcing on malleable systems. Since malleable systems have a wider range of dynamical states available by changing parameters, a time-dependent change in the parameters, induced by the noise or forcing, can lead to transitions between several different states. The complex and sensitive dependence on parameters would mean that even small amplitude changes could lead to drastic fluctuations. For the Watts-Strogatz networks, multistability can complicate the dependence on external inputs, and make the effects dependent on the timing of perturbations, as different states, all of which coexist, can react differently to the parameter changes. Understanding these behaviors is important, for instance, in the context of neural systems, where external influences are common and where temporal fluctuations are essential.

Future research is also needed to fully describe the mechanism for the STS fluctuations. An interesting possibility could be to extend the synchrony alignment function \cite{skardal2014optimal} to weakly synchronized regimes. Another promising approach would also be to apply the formalism introduced in \cite{muller2021algebraic,budzinski2022geometry}. This would be an important theoretical contribution for the understanding of phase synchronization in oscillator networks and for the role of each unit in a network.

To summarize, the increased magnitude and complexity of dynamical malleability is a general phenomenon in finite-size systems and can be expected to occur in real-world systems. 

\begin{acknowledgments}
We would like to thank Jan Freund and Arkady Pikovsky for helpful discussions. K.L.R. was supported by the German Academic Exchange Service (DAAD). R.C.B. and L.E.M. acknowledge the support by BrainsCAN at Western University through the Canada First Research Excellence Fund (CFREF), the NSF through a NeuroNex award (\#2015276), SPIRITS 2020 of Kyoto University, Compute Ontario (computeontario.ca), Compute Canada (computecanada.ca), and the Western Academy for Advanced Research. R.C.B gratefully acknowledges the Western Institute for Neuroscience Clinical Research Postdoctoral Fellowship. B. R. R. B. acknowledges the financial support of the São Paulo Research Foundation (FAPESP, Brazil) Grants Nos. 2018/03211-6 and 2021/09839-0. The simulations were performed at the HPC Cluster CARL, located at the University of Oldenburg (Germany) and funded by the DFG through its Major Research Instrumentation Program (INST 184/157-1 FUGG) and the Ministry of Science and Culture (MWK) of the Lower Saxony State, Germany.
\end{acknowledgments}


\begin{thebibliography}{67}%
\makeatletter
\providecommand \@ifxundefined [1]{%
 \@ifx{#1\undefined}
}%
\providecommand \@ifnum [1]{%
 \ifnum #1\expandafter \@firstoftwo
 \else \expandafter \@secondoftwo
 \fi
}%
\providecommand \@ifx [1]{%
 \ifx #1\expandafter \@firstoftwo
 \else \expandafter \@secondoftwo
 \fi
}%
\providecommand \natexlab [1]{#1}%
\providecommand \enquote  [1]{``#1''}%
\providecommand \bibnamefont  [1]{#1}%
\providecommand \bibfnamefont [1]{#1}%
\providecommand \citenamefont [1]{#1}%
\providecommand \href@noop [0]{\@secondoftwo}%
\providecommand \href [0]{\begingroup \@sanitize@url \@href}%
\providecommand \@href[1]{\@@startlink{#1}\@@href}%
\providecommand \@@href[1]{\endgroup#1\@@endlink}%
\providecommand \@sanitize@url [0]{\catcode `\\12\catcode `\$12\catcode
  `\&12\catcode `\#12\catcode `\^12\catcode `\_12\catcode `\%12\relax}%
\providecommand \@@startlink[1]{}%
\providecommand \@@endlink[0]{}%
\providecommand \url  [0]{\begingroup\@sanitize@url \@url }%
\providecommand \@url [1]{\endgroup\@href {#1}{\urlprefix }}%
\providecommand \urlprefix  [0]{URL }%
\providecommand \Eprint [0]{\href }%
\providecommand \doibase [0]{https://doi.org/}%
\providecommand \selectlanguage [0]{\@gobble}%
\providecommand \bibinfo  [0]{\@secondoftwo}%
\providecommand \bibfield  [0]{\@secondoftwo}%
\providecommand \translation [1]{[#1]}%
\providecommand \BibitemOpen [0]{}%
\providecommand \bibitemStop [0]{}%
\providecommand \bibitemNoStop [0]{.\EOS\space}%
\providecommand \EOS [0]{\spacefactor3000\relax}%
\providecommand \BibitemShut  [1]{\csname bibitem#1\endcsname}%
\let\auto@bib@innerbib\@empty
\bibitem [{\citenamefont {Motter}\ \emph {et~al.}(2013)\citenamefont {Motter},
  \citenamefont {Myers}, \citenamefont {Anghel},\ and\ \citenamefont
  {Nishikawa}}]{motter2013spontaneous}%
  \BibitemOpen
  \bibfield  {author} {\bibinfo {author} {\bibfnamefont {A.~E.}\ \bibnamefont
  {Motter}}, \bibinfo {author} {\bibfnamefont {S.~A.}\ \bibnamefont {Myers}},
  \bibinfo {author} {\bibfnamefont {M.}~\bibnamefont {Anghel}},\ and\ \bibinfo
  {author} {\bibfnamefont {T.}~\bibnamefont {Nishikawa}},\ }\bibfield  {title}
  {\bibinfo {title} {{Spontaneous synchrony in power-grid networks}},\ }\href
  {https://doi.org/10.1038/nphys2535} {\bibfield  {journal} {\bibinfo
  {journal} {Nature Physics}\ }\textbf {\bibinfo {volume} {9}},\ \bibinfo
  {pages} {191} (\bibinfo {year} {2013})}\BibitemShut {NoStop}%
\bibitem [{\citenamefont {Dunne}\ \emph {et~al.}(2002)\citenamefont {Dunne},
  \citenamefont {Williams},\ and\ \citenamefont {Martinez}}]{dunne2002foodweb}%
  \BibitemOpen
  \bibfield  {author} {\bibinfo {author} {\bibfnamefont {J.~A.}\ \bibnamefont
  {Dunne}}, \bibinfo {author} {\bibfnamefont {R.~J.}\ \bibnamefont
  {Williams}},\ and\ \bibinfo {author} {\bibfnamefont {N.~D.}\ \bibnamefont
  {Martinez}},\ }\bibfield  {title} {\bibinfo {title} {{Food-web structure and
  network theory: The role of connectance and size}},\ }\href
  {https://doi.org/10.1073/pnas.192407699} {\bibfield  {journal} {\bibinfo
  {journal} {Proceedings of the National Academy of Sciences}\ }\textbf
  {\bibinfo {volume} {99}},\ \bibinfo {pages} {12917} (\bibinfo {year}
  {2002})}\BibitemShut {NoStop}%
\bibitem [{\citenamefont {Crotty}\ \emph {et~al.}(2010)\citenamefont {Crotty},
  \citenamefont {Schult},\ and\ \citenamefont {Segall}}]{crotty2010josephson}%
  \BibitemOpen
  \bibfield  {author} {\bibinfo {author} {\bibfnamefont {P.}~\bibnamefont
  {Crotty}}, \bibinfo {author} {\bibfnamefont {D.}~\bibnamefont {Schult}},\
  and\ \bibinfo {author} {\bibfnamefont {K.}~\bibnamefont {Segall}},\
  }\bibfield  {title} {\bibinfo {title} {{Josephson junction simulation of
  neurons}},\ }\href {https://doi.org/10.1103/physreve.82.011914} {\bibfield
  {journal} {\bibinfo  {journal} {Physical Review E}\ }\textbf {\bibinfo
  {volume} {82}},\ \bibinfo {pages} {011914} (\bibinfo {year}
  {2010})}\BibitemShut {NoStop}%
\bibitem [{\citenamefont {Nixon}\ \emph {et~al.}(2011)\citenamefont {Nixon},
  \citenamefont {Friedman}, \citenamefont {Ronen}, \citenamefont {Friesem},
  \citenamefont {Davidson},\ and\ \citenamefont
  {Kanter}}]{nixon2011synchronized}%
  \BibitemOpen
  \bibfield  {author} {\bibinfo {author} {\bibfnamefont {M.}~\bibnamefont
  {Nixon}}, \bibinfo {author} {\bibfnamefont {M.}~\bibnamefont {Friedman}},
  \bibinfo {author} {\bibfnamefont {E.}~\bibnamefont {Ronen}}, \bibinfo
  {author} {\bibfnamefont {A.~A.}\ \bibnamefont {Friesem}}, \bibinfo {author}
  {\bibfnamefont {N.}~\bibnamefont {Davidson}},\ and\ \bibinfo {author}
  {\bibfnamefont {I.}~\bibnamefont {Kanter}},\ }\bibfield  {title} {\bibinfo
  {title} {{Synchronized Cluster Formation in Coupled Laser Networks}},\ }\href
  {https://doi.org/10.1103/physrevlett.106.223901} {\bibfield  {journal}
  {\bibinfo  {journal} {Physical Review Letters}\ }\textbf {\bibinfo {volume}
  {106}},\ \bibinfo {pages} {223901} (\bibinfo {year} {2011})}\BibitemShut
  {NoStop}%
\bibitem [{\citenamefont {Varela}\ \emph {et~al.}(2001)\citenamefont {Varela},
  \citenamefont {Lachaux}, \citenamefont {Rodriguez},\ and\ \citenamefont
  {Martinerie}}]{varela2001brainweb}%
  \BibitemOpen
  \bibfield  {author} {\bibinfo {author} {\bibfnamefont {F.}~\bibnamefont
  {Varela}}, \bibinfo {author} {\bibfnamefont {J.~P.}\ \bibnamefont {Lachaux}},
  \bibinfo {author} {\bibfnamefont {E.}~\bibnamefont {Rodriguez}},\ and\
  \bibinfo {author} {\bibfnamefont {J.}~\bibnamefont {Martinerie}},\ }\bibfield
   {title} {\bibinfo {title} {{The brainweb: phase synchronization and
  large-scale integration.}},\ }\href {https://doi.org/10.1038/35067550}
  {\bibfield  {journal} {\bibinfo  {journal} {Nature Reviews. Neuroscience}\
  }\textbf {\bibinfo {volume} {2}},\ \bibinfo {pages} {229} (\bibinfo {year}
  {2001})}\BibitemShut {NoStop}%
\bibitem [{\citenamefont {Witthaut}\ \emph {et~al.}(2022)\citenamefont
  {Witthaut}, \citenamefont {Hellmann}, \citenamefont {Kurths}, \citenamefont
  {Kettemann}, \citenamefont {Meyer-Ortmanns},\ and\ \citenamefont
  {Timme}}]{witthaut2022collective}%
  \BibitemOpen
  \bibfield  {author} {\bibinfo {author} {\bibfnamefont {D.}~\bibnamefont
  {Witthaut}}, \bibinfo {author} {\bibfnamefont {F.}~\bibnamefont {Hellmann}},
  \bibinfo {author} {\bibfnamefont {J.}~\bibnamefont {Kurths}}, \bibinfo
  {author} {\bibfnamefont {S.}~\bibnamefont {Kettemann}}, \bibinfo {author}
  {\bibfnamefont {H.}~\bibnamefont {Meyer-Ortmanns}},\ and\ \bibinfo {author}
  {\bibfnamefont {M.}~\bibnamefont {Timme}},\ }\bibfield  {title} {\bibinfo
  {title} {{Collective nonlinear dynamics and self-organization in
  decentralized power grids}},\ }\href
  {https://doi.org/10.1103/revmodphys.94.015005} {\bibfield  {journal}
  {\bibinfo  {journal} {Reviews of Modern Physics}\ }\textbf {\bibinfo {volume}
  {94}},\ \bibinfo {pages} {015005} (\bibinfo {year} {2022})}\BibitemShut
  {NoStop}%
\bibitem [{\citenamefont {Fries}(2015)}]{fries2015rhythms}%
  \BibitemOpen
  \bibfield  {author} {\bibinfo {author} {\bibfnamefont {P.}~\bibnamefont
  {Fries}},\ }\bibfield  {title} {\bibinfo {title} {{Rhythms for Cognition:
  Communication through Coherence.}},\ }\href
  {https://doi.org/10.1016/j.neuron.2015.09.034} {\bibfield  {journal}
  {\bibinfo  {journal} {Neuron}\ }\textbf {\bibinfo {volume} {88}},\ \bibinfo
  {pages} {220} (\bibinfo {year} {2015})}\BibitemShut {NoStop}%
\bibitem [{\citenamefont {Singer}(1999)}]{singer1999neuronal}%
  \BibitemOpen
  \bibfield  {author} {\bibinfo {author} {\bibfnamefont {W.}~\bibnamefont
  {Singer}},\ }\bibfield  {title} {\bibinfo {title} {{Neuronal synchrony: a
  versatile code for the definition of relations?}},\ }\href
  {https://doi.org/10.1016/s0896-6273(00)80821-1} {\bibfield  {journal}
  {\bibinfo  {journal} {Neuron}\ }\textbf {\bibinfo {volume} {24}},\ \bibinfo
  {pages} {49} (\bibinfo {year} {1999})}\BibitemShut {NoStop}%
\bibitem [{\citenamefont {Pikovsky}\ \emph {et~al.}(2001)\citenamefont
  {Pikovsky}, \citenamefont {Rosenblum},\ and\ \citenamefont
  {Kurths}}]{pikovsky2001synchronization}%
  \BibitemOpen
  \bibfield  {author} {\bibinfo {author} {\bibfnamefont {A.}~\bibnamefont
  {Pikovsky}}, \bibinfo {author} {\bibfnamefont {M.}~\bibnamefont
  {Rosenblum}},\ and\ \bibinfo {author} {\bibfnamefont {J.}~\bibnamefont
  {Kurths}},\ }\href {https://doi.org/10.1119/1.1475332} {\emph {\bibinfo
  {title} {{Synchronization: A Universal Concept in Nonlinear Science}}}},\
  Vol.~\bibinfo {volume} {70}\ (\bibinfo  {publisher} {Cambridge University
  Press},\ \bibinfo {year} {2001})\BibitemShut {NoStop}%
\bibitem [{\citenamefont {Arenas}\ \emph {et~al.}(2008)\citenamefont {Arenas},
  \citenamefont {Díaz-Guilera}, \citenamefont {Kurths}, \citenamefont
  {Moreno},\ and\ \citenamefont {Zhou}}]{arenas2008synchronization}%
  \BibitemOpen
  \bibfield  {author} {\bibinfo {author} {\bibfnamefont {A.}~\bibnamefont
  {Arenas}}, \bibinfo {author} {\bibfnamefont {A.}~\bibnamefont
  {Díaz-Guilera}}, \bibinfo {author} {\bibfnamefont {J.}~\bibnamefont
  {Kurths}}, \bibinfo {author} {\bibfnamefont {Y.}~\bibnamefont {Moreno}},\
  and\ \bibinfo {author} {\bibfnamefont {C.}~\bibnamefont {Zhou}},\ }\bibfield
  {title} {\bibinfo {title} {{Synchronization in complex networks}},\ }\href
  {https://doi.org/10.1016/j.physrep.2008.09.002} {\bibfield  {journal}
  {\bibinfo  {journal} {Physics Reports}\ }\textbf {\bibinfo {volume} {469}},\
  \bibinfo {pages} {93} (\bibinfo {year} {2008})},\ \Eprint
  {https://arxiv.org/abs/0805.2976} {0805.2976} \BibitemShut {NoStop}%
\bibitem [{\citenamefont {Mitra}\ \emph {et~al.}(2017)\citenamefont {Mitra},
  \citenamefont {Choudhary}, \citenamefont {Sinha}, \citenamefont {Kurths},\
  and\ \citenamefont {Donner}}]{mitra2017multiple}%
  \BibitemOpen
  \bibfield  {author} {\bibinfo {author} {\bibfnamefont {C.}~\bibnamefont
  {Mitra}}, \bibinfo {author} {\bibfnamefont {A.}~\bibnamefont {Choudhary}},
  \bibinfo {author} {\bibfnamefont {S.}~\bibnamefont {Sinha}}, \bibinfo
  {author} {\bibfnamefont {J.}~\bibnamefont {Kurths}},\ and\ \bibinfo {author}
  {\bibfnamefont {R.~V.}\ \bibnamefont {Donner}},\ }\bibfield  {title}
  {\bibinfo {title} {{Multiple-node basin stability in complex dynamical
  networks}},\ }\href {https://doi.org/10.1103/physreve.95.032317} {\bibfield
  {journal} {\bibinfo  {journal} {Physical Review E}\ }\textbf {\bibinfo
  {volume} {95}},\ \bibinfo {pages} {032317} (\bibinfo {year} {2017})},\
  \Eprint {https://arxiv.org/abs/1612.06015} {1612.06015} \BibitemShut
  {NoStop}%
\bibitem [{\citenamefont {Medeiros}\ \emph {et~al.}(2018)\citenamefont
  {Medeiros}, \citenamefont {Medrano-T}, \citenamefont {Caldas},\ and\
  \citenamefont {Feudel}}]{medeiros2018boundaries}%
  \BibitemOpen
  \bibfield  {author} {\bibinfo {author} {\bibfnamefont {E.~S.}\ \bibnamefont
  {Medeiros}}, \bibinfo {author} {\bibfnamefont {R.~O.}\ \bibnamefont
  {Medrano-T}}, \bibinfo {author} {\bibfnamefont {I.~L.}\ \bibnamefont
  {Caldas}},\ and\ \bibinfo {author} {\bibfnamefont {U.}~\bibnamefont
  {Feudel}},\ }\bibfield  {title} {\bibinfo {title} {{Boundaries of
  synchronization in oscillator networks}},\ }\href
  {https://doi.org/10.1103/physreve.98.030201} {\bibfield  {journal} {\bibinfo
  {journal} {Physical Review E}\ }\textbf {\bibinfo {volume} {98}},\ \bibinfo
  {pages} {030201} (\bibinfo {year} {2018})}\BibitemShut {NoStop}%
\bibitem [{\citenamefont {Witthaut}\ and\ \citenamefont
  {Timme}(2012)}]{witthaut2012braess}%
  \BibitemOpen
  \bibfield  {author} {\bibinfo {author} {\bibfnamefont {D.}~\bibnamefont
  {Witthaut}}\ and\ \bibinfo {author} {\bibfnamefont {M.}~\bibnamefont
  {Timme}},\ }\bibfield  {title} {\bibinfo {title} {{Braess's paradox in
  oscillator networks, desynchronization and power outage}},\ }\href
  {https://doi.org/10.1088/1367-2630/14/8/083036} {\bibfield  {journal}
  {\bibinfo  {journal} {New Journal of Physics}\ }\textbf {\bibinfo {volume}
  {14}},\ \bibinfo {pages} {083036} (\bibinfo {year} {2012})}\BibitemShut
  {NoStop}%
\bibitem [{\citenamefont {Coletta}\ and\ \citenamefont
  {Jacquod}(2016)}]{coletta2016linear}%
  \BibitemOpen
  \bibfield  {author} {\bibinfo {author} {\bibfnamefont {T.}~\bibnamefont
  {Coletta}}\ and\ \bibinfo {author} {\bibfnamefont {P.}~\bibnamefont
  {Jacquod}},\ }\bibfield  {title} {\bibinfo {title} {{Linear stability and the
  Braess paradox in coupled-oscillator networks and electric power grids}},\
  }\href {https://doi.org/10.1103/physreve.93.032222} {\bibfield  {journal}
  {\bibinfo  {journal} {Physical Review E}\ }\textbf {\bibinfo {volume} {93}},\
  \bibinfo {pages} {032222} (\bibinfo {year} {2016})}\BibitemShut {NoStop}%
\bibitem [{\citenamefont {Mihara}\ \emph {et~al.}(2022)\citenamefont {Mihara},
  \citenamefont {Medeiros}, \citenamefont {Zakharova},\ and\ \citenamefont
  {Medrano-T}}]{mihara2022sparsity}%
  \BibitemOpen
  \bibfield  {author} {\bibinfo {author} {\bibfnamefont {A.}~\bibnamefont
  {Mihara}}, \bibinfo {author} {\bibfnamefont {E.~S.}\ \bibnamefont
  {Medeiros}}, \bibinfo {author} {\bibfnamefont {A.}~\bibnamefont
  {Zakharova}},\ and\ \bibinfo {author} {\bibfnamefont {R.~O.}\ \bibnamefont
  {Medrano-T}},\ }\bibfield  {title} {\bibinfo {title} {{Sparsity-driven
  synchronization in oscillator networks}},\ }\href
  {https://doi.org/10.1063/5.0074008} {\bibfield  {journal} {\bibinfo
  {journal} {Chaos: An Interdisciplinary Journal of Nonlinear Science}\
  }\textbf {\bibinfo {volume} {32}},\ \bibinfo {pages} {033114} (\bibinfo
  {year} {2022})},\ \Eprint {https://arxiv.org/abs/2111.13583} {2111.13583}
  \BibitemShut {NoStop}%
\bibitem [{\citenamefont {Menck}\ \emph {et~al.}(2014)\citenamefont {Menck},
  \citenamefont {Heitzig}, \citenamefont {Kurths},\ and\ \citenamefont
  {Schellnhuber}}]{menck2014how}%
  \BibitemOpen
  \bibfield  {author} {\bibinfo {author} {\bibfnamefont {P.~J.}\ \bibnamefont
  {Menck}}, \bibinfo {author} {\bibfnamefont {J.}~\bibnamefont {Heitzig}},
  \bibinfo {author} {\bibfnamefont {J.}~\bibnamefont {Kurths}},\ and\ \bibinfo
  {author} {\bibfnamefont {H.~J.}\ \bibnamefont {Schellnhuber}},\ }\bibfield
  {title} {\bibinfo {title} {{How dead ends undermine power grid stability}},\
  }\href {https://doi.org/10.1038/ncomms4969} {\bibfield  {journal} {\bibinfo
  {journal} {Nature Communications}\ }\textbf {\bibinfo {volume} {5}},\
  \bibinfo {pages} {3969} (\bibinfo {year} {2014})}\BibitemShut {NoStop}%
\bibitem [{\citenamefont {Halekotte}\ and\ \citenamefont
  {Feudel}(2020)}]{halekotte2020minimal}%
  \BibitemOpen
  \bibfield  {author} {\bibinfo {author} {\bibfnamefont {L.}~\bibnamefont
  {Halekotte}}\ and\ \bibinfo {author} {\bibfnamefont {U.}~\bibnamefont
  {Feudel}},\ }\bibfield  {title} {\bibinfo {title} {{Minimal fatal shocks in
  multistable complex networks}},\ }\href
  {https://doi.org/10.1038/s41598-020-68805-6} {\bibfield  {journal} {\bibinfo
  {journal} {Scientific Reports}\ }\textbf {\bibinfo {volume} {10}},\ \bibinfo
  {pages} {11783} (\bibinfo {year} {2020})}\BibitemShut {NoStop}%
\bibitem [{\citenamefont {Manik}\ \emph {et~al.}(2017)\citenamefont {Manik},
  \citenamefont {Rohden}, \citenamefont {Ronellenfitsch}, \citenamefont
  {Zhang}, \citenamefont {Hallerberg}, \citenamefont {Witthaut},\ and\
  \citenamefont {Timme}}]{manik2017network}%
  \BibitemOpen
  \bibfield  {author} {\bibinfo {author} {\bibfnamefont {D.}~\bibnamefont
  {Manik}}, \bibinfo {author} {\bibfnamefont {M.}~\bibnamefont {Rohden}},
  \bibinfo {author} {\bibfnamefont {H.}~\bibnamefont {Ronellenfitsch}},
  \bibinfo {author} {\bibfnamefont {X.}~\bibnamefont {Zhang}}, \bibinfo
  {author} {\bibfnamefont {S.}~\bibnamefont {Hallerberg}}, \bibinfo {author}
  {\bibfnamefont {D.}~\bibnamefont {Witthaut}},\ and\ \bibinfo {author}
  {\bibfnamefont {M.}~\bibnamefont {Timme}},\ }\bibfield  {title} {\bibinfo
  {title} {{Network susceptibilities: Theory and applications}},\ }\href
  {https://doi.org/10.1103/physreve.95.012319} {\bibfield  {journal} {\bibinfo
  {journal} {Physical Review E}\ }\textbf {\bibinfo {volume} {95}},\ \bibinfo
  {pages} {012319} (\bibinfo {year} {2017})},\ \Eprint
  {https://arxiv.org/abs/1609.04310} {1609.04310} \BibitemShut {NoStop}%
\bibitem [{\citenamefont {Buenda}\ \emph {et~al.}(2022)\citenamefont {Buenda},
  \citenamefont {Villegas}, \citenamefont {Burioni},\ and\ \citenamefont
  {Muoz}}]{buendia2022the}%
  \BibitemOpen
  \bibfield  {author} {\bibinfo {author} {\bibfnamefont {V.}~\bibnamefont
  {Buenda}}, \bibinfo {author} {\bibfnamefont {P.}~\bibnamefont {Villegas}},
  \bibinfo {author} {\bibfnamefont {R.}~\bibnamefont {Burioni}},\ and\ \bibinfo
  {author} {\bibfnamefont {M.~A.}\ \bibnamefont {Muoz}},\ }\bibfield  {title}
  {\bibinfo {title} {{The broad edge of synchronization: Griffiths effects and
  collective phenomena in brain networks}},\ }\href
  {https://doi.org/10.1098/rsta.2020.0424} {\bibfield  {journal} {\bibinfo
  {journal} {Philosophical Transactions of the Royal Society A}\ }\textbf
  {\bibinfo {volume} {380}},\ \bibinfo {pages} {20200424} (\bibinfo {year}
  {2022})}\BibitemShut {NoStop}%
\bibitem [{\citenamefont {Kuramoto}(1975)}]{kuramoto1975self}%
  \BibitemOpen
  \bibfield  {author} {\bibinfo {author} {\bibfnamefont {Y.}~\bibnamefont
  {Kuramoto}},\ }\bibfield  {title} {\bibinfo {title} {{Self-entrainment of a
  population of coupled non-linear oscillators}},\ }\href
  {https://doi.org/10.1007/bfb0013365} {\bibfield  {journal} {\bibinfo
  {journal} {International Symposium on Mathematical Problems in Theoretical
  Physics. Lecture Notes in Physics}\ }\textbf {\bibinfo {volume} {39}},\
  \bibinfo {pages} {420} (\bibinfo {year} {1975})}\BibitemShut {NoStop}%
\bibitem [{\citenamefont {Acebrón}\ \emph {et~al.}(2005)\citenamefont
  {Acebrón}, \citenamefont {Bonilla}, \citenamefont {Vicente}, \citenamefont
  {Ritort},\ and\ \citenamefont {Spigler}}]{acebron2005kuramoto}%
  \BibitemOpen
  \bibfield  {author} {\bibinfo {author} {\bibfnamefont {J.~A.}\ \bibnamefont
  {Acebrón}}, \bibinfo {author} {\bibfnamefont {L.~L.}\ \bibnamefont
  {Bonilla}}, \bibinfo {author} {\bibfnamefont {C.~J.}\ \bibnamefont
  {Vicente}}, \bibinfo {author} {\bibfnamefont {F.}~\bibnamefont {Ritort}},\
  and\ \bibinfo {author} {\bibfnamefont {R.}~\bibnamefont {Spigler}},\
  }\bibfield  {title} {\bibinfo {title} {{The Kuramoto model: A simple paradigm
  for synchronization phenomena}},\ }\href
  {https://doi.org/10.1103/revmodphys.77.137} {\bibfield  {journal} {\bibinfo
  {journal} {Reviews of Modern Physics}\ }\textbf {\bibinfo {volume} {77}},\
  \bibinfo {pages} {137} (\bibinfo {year} {2005})}\BibitemShut {NoStop}%
\bibitem [{\citenamefont {Rodrigues}\ \emph {et~al.}(2016)\citenamefont
  {Rodrigues}, \citenamefont {Peron}, \citenamefont {Ji},\ and\ \citenamefont
  {Kurths}}]{rodrigues2016the}%
  \BibitemOpen
  \bibfield  {author} {\bibinfo {author} {\bibfnamefont {F.~A.}\ \bibnamefont
  {Rodrigues}}, \bibinfo {author} {\bibfnamefont {T.~K.}\ \bibnamefont
  {Peron}}, \bibinfo {author} {\bibfnamefont {P.}~\bibnamefont {Ji}},\ and\
  \bibinfo {author} {\bibfnamefont {J.}~\bibnamefont {Kurths}},\ }\bibfield
  {title} {\bibinfo {title} {{The Kuramoto model in complex networks}},\ }\href
  {https://doi.org/10.1016/j.physrep.2015.10.008} {\bibfield  {journal}
  {\bibinfo  {journal} {Physics Reports}\ }\textbf {\bibinfo {volume} {610}},\
  \bibinfo {pages} {1} (\bibinfo {year} {2016})}\BibitemShut {NoStop}%
\bibitem [{\citenamefont {Ponce-Alvarez}\ \emph {et~al.}(2015)\citenamefont
  {Ponce-Alvarez}, \citenamefont {Deco}, \citenamefont {Hagmann}, \citenamefont
  {Romani}, \citenamefont {Mantini},\ and\ \citenamefont
  {Corbetta}}]{ponce-alvarez2015restingstate}%
  \BibitemOpen
  \bibfield  {author} {\bibinfo {author} {\bibfnamefont {A.}~\bibnamefont
  {Ponce-Alvarez}}, \bibinfo {author} {\bibfnamefont {G.}~\bibnamefont {Deco}},
  \bibinfo {author} {\bibfnamefont {P.}~\bibnamefont {Hagmann}}, \bibinfo
  {author} {\bibfnamefont {G.~L.}\ \bibnamefont {Romani}}, \bibinfo {author}
  {\bibfnamefont {D.}~\bibnamefont {Mantini}},\ and\ \bibinfo {author}
  {\bibfnamefont {M.}~\bibnamefont {Corbetta}},\ }\bibfield  {title} {\bibinfo
  {title} {{Resting-State Temporal Synchronization Networks Emerge from
  Connectivity Topology and Heterogeneity}},\ }\href
  {https://doi.org/10.1371/journal.pcbi.1004100} {\bibfield  {journal}
  {\bibinfo  {journal} {PLOS Computational Biology}\ }\textbf {\bibinfo
  {volume} {11}},\ \bibinfo {pages} {e1004100} (\bibinfo {year}
  {2015})}\BibitemShut {NoStop}%
\bibitem [{\citenamefont {Cabral}\ \emph {et~al.}(2011)\citenamefont {Cabral},
  \citenamefont {Hugues}, \citenamefont {Sporns},\ and\ \citenamefont
  {Deco}}]{cabral2011role}%
  \BibitemOpen
  \bibfield  {author} {\bibinfo {author} {\bibfnamefont {J.}~\bibnamefont
  {Cabral}}, \bibinfo {author} {\bibfnamefont {E.}~\bibnamefont {Hugues}},
  \bibinfo {author} {\bibfnamefont {O.}~\bibnamefont {Sporns}},\ and\ \bibinfo
  {author} {\bibfnamefont {G.}~\bibnamefont {Deco}},\ }\bibfield  {title}
  {\bibinfo {title} {{Role of local network oscillations in resting-state
  functional connectivity.}},\ }\href
  {https://doi.org/10.1016/j.neuroimage.2011.04.010} {\bibfield  {journal}
  {\bibinfo  {journal} {Neuroimage}\ }\textbf {\bibinfo {volume} {57}},\
  \bibinfo {pages} {130} (\bibinfo {year} {2011})}\BibitemShut {NoStop}%
\bibitem [{\citenamefont {Josephson}(1964)}]{josephson1964coupled}%
  \BibitemOpen
  \bibfield  {author} {\bibinfo {author} {\bibfnamefont {B.~D.}\ \bibnamefont
  {Josephson}},\ }\bibfield  {title} {\bibinfo {title} {{Coupled
  Superconductors}},\ }\href {https://doi.org/10.1103/revmodphys.36.216}
  {\bibfield  {journal} {\bibinfo  {journal} {Reviews of Modern Physics}\
  }\textbf {\bibinfo {volume} {36}},\ \bibinfo {pages} {216} (\bibinfo {year}
  {1964})}\BibitemShut {NoStop}%
\bibitem [{\citenamefont {Marek}\ and\ \citenamefont
  {Stuchl}(1975)}]{marek1975synchronization}%
  \BibitemOpen
  \bibfield  {author} {\bibinfo {author} {\bibfnamefont {M.}~\bibnamefont
  {Marek}}\ and\ \bibinfo {author} {\bibfnamefont {I.}~\bibnamefont {Stuchl}},\
  }\bibfield  {title} {\bibinfo {title} {{Synchronization in two interacting
  oscillatory systems}},\ }\href {https://doi.org/10.1016/0301-4622(75)80016-0}
  {\bibfield  {journal} {\bibinfo  {journal} {Biophysical Chemistry}\ }\textbf
  {\bibinfo {volume} {3}},\ \bibinfo {pages} {241} (\bibinfo {year}
  {1975})}\BibitemShut {NoStop}%
\bibitem [{\citenamefont {Neu}(1979)}]{neu1979chemical}%
  \BibitemOpen
  \bibfield  {author} {\bibinfo {author} {\bibfnamefont {J.~C.}\ \bibnamefont
  {Neu}},\ }\bibfield  {title} {\bibinfo {title} {{Chemical Waves and the
  Diffusive Coupling of Limit Cycle Oscillators}},\ }\href
  {https://doi.org/10.1137/0136038} {\bibfield  {journal} {\bibinfo  {journal}
  {SIAM Journal on Applied Mathematics}\ }\textbf {\bibinfo {volume} {36}},\
  \bibinfo {pages} {509} (\bibinfo {year} {1979})}\BibitemShut {NoStop}%
\bibitem [{\citenamefont {Kuramoto}(1984)}]{kuramoto1984chemical}%
  \BibitemOpen
  \bibfield  {author} {\bibinfo {author} {\bibfnamefont {Y.}~\bibnamefont
  {Kuramoto}},\ }\href {https://doi.org/10.1007/978-3-642-69689-3} {\emph
  {\bibinfo {title} {{Chemical Oscillations, Waves, and Turbulence}}}},\
  \bibinfo {series} {Springer Series in Synergetics}, Vol.~\bibinfo {volume}
  {19}\ (\bibinfo  {publisher} {Springer Berlin Heidelberg},\ \bibinfo
  {address} {Berlin, Heidelberg},\ \bibinfo {year} {1984})\BibitemShut
  {NoStop}%
\bibitem [{\citenamefont {Strogatz}(2000)}]{strogatz2000from}%
  \BibitemOpen
  \bibfield  {author} {\bibinfo {author} {\bibfnamefont {S.~H.}\ \bibnamefont
  {Strogatz}},\ }\bibfield  {title} {\bibinfo {title} {{From Kuramoto to
  Crawford: exploring the onset of synchronization in populations of coupled
  oscillators}},\ }\href {https://doi.org/10.1016/s0167-2789(00)00094-4}
  {\bibfield  {journal} {\bibinfo  {journal} {Physica D: Nonlinear Phenomena}\
  }\textbf {\bibinfo {volume} {143}},\ \bibinfo {pages} {1} (\bibinfo {year}
  {2000})}\BibitemShut {NoStop}%
\bibitem [{\citenamefont {Kuramoto}\ and\ \citenamefont
  {Nakao}(2019)}]{kuramoto2019on}%
  \BibitemOpen
  \bibfield  {author} {\bibinfo {author} {\bibfnamefont {Y.}~\bibnamefont
  {Kuramoto}}\ and\ \bibinfo {author} {\bibfnamefont {H.}~\bibnamefont
  {Nakao}},\ }\bibfield  {title} {\bibinfo {title} {{On the concept of
  dynamical reduction: the case of coupled oscillators}},\ }\href
  {https://doi.org/10.1098/rsta.2019.0041} {\bibfield  {journal} {\bibinfo
  {journal} {Philosophical Transactions of the Royal Society A}\ }\textbf
  {\bibinfo {volume} {377}},\ \bibinfo {pages} {20190041} (\bibinfo {year}
  {2019})}\BibitemShut {NoStop}%
\bibitem [{\citenamefont {Watts}\ and\ \citenamefont
  {Strogatz}(1998)}]{watts1998collective}%
  \BibitemOpen
  \bibfield  {author} {\bibinfo {author} {\bibfnamefont {D.~J.}\ \bibnamefont
  {Watts}}\ and\ \bibinfo {author} {\bibfnamefont {S.~H.}\ \bibnamefont
  {Strogatz}},\ }\bibfield  {title} {\bibinfo {title} {{Collective dynamics of
  'small-world' networks}},\ }\href {https://doi.org/10.1038/30918} {\bibfield
  {journal} {\bibinfo  {journal} {Nature}\ }\textbf {\bibinfo {volume} {393}},\
  \bibinfo {pages} {440} (\bibinfo {year} {1998})}\BibitemShut {NoStop}%
\bibitem [{\citenamefont {Rogers}\ and\ \citenamefont
  {Wille}(1996)}]{rogers1996phasetransitions}%
  \BibitemOpen
  \bibfield  {author} {\bibinfo {author} {\bibfnamefont {J.~L.}\ \bibnamefont
  {Rogers}}\ and\ \bibinfo {author} {\bibfnamefont {L.~T.}\ \bibnamefont
  {Wille}},\ }\bibfield  {title} {\bibinfo {title} {{Phase transitions in
  nonlinear oscillator chains}},\ }\href
  {https://doi.org/10.1103/physreve.54.r2193} {\bibfield  {journal} {\bibinfo
  {journal} {Physical Review E}\ }\textbf {\bibinfo {volume} {54}},\ \bibinfo
  {pages} {R2193} (\bibinfo {year} {1996})}\BibitemShut {NoStop}%
\bibitem [{\citenamefont {Skardal}\ and\ \citenamefont
  {Arenas}(2020)}]{skardal2020higher}%
  \BibitemOpen
  \bibfield  {author} {\bibinfo {author} {\bibfnamefont {P.~S.}\ \bibnamefont
  {Skardal}}\ and\ \bibinfo {author} {\bibfnamefont {A.}~\bibnamefont
  {Arenas}},\ }\bibfield  {title} {\bibinfo {title} {{Higher order interactions
  in complex networks of phase oscillators promote abrupt synchronization
  switching}},\ }\href {https://doi.org/10.1038/s42005-020-00485-0} {\bibfield
  {journal} {\bibinfo  {journal} {Communications Physics}\ }\textbf {\bibinfo
  {volume} {3}},\ \bibinfo {pages} {218} (\bibinfo {year} {2020})}\BibitemShut
  {NoStop}%
\bibitem [{\citenamefont {Hong}\ \emph {et~al.}(2013)\citenamefont {Hong},
  \citenamefont {Um},\ and\ \citenamefont {Park}}]{hong2013link}%
  \BibitemOpen
  \bibfield  {author} {\bibinfo {author} {\bibfnamefont {H.}~\bibnamefont
  {Hong}}, \bibinfo {author} {\bibfnamefont {J.}~\bibnamefont {Um}},\ and\
  \bibinfo {author} {\bibfnamefont {H.}~\bibnamefont {Park}},\ }\bibfield
  {title} {\bibinfo {title} {{Link-disorder fluctuation effects on
  synchronization in random networks}},\ }\href
  {https://doi.org/10.1103/physreve.87.042105} {\bibfield  {journal} {\bibinfo
  {journal} {Physical Review E - Statistical, Nonlinear, and Soft Matter
  Physics}\ }\textbf {\bibinfo {volume} {87}},\ \bibinfo {pages} {1} (\bibinfo
  {year} {2013})},\ \Eprint {https://arxiv.org/abs/1304.0610} {1304.0610}
  \BibitemShut {NoStop}%
\bibitem [{\citenamefont {Hong}\ \emph {et~al.}(2002)\citenamefont {Hong},
  \citenamefont {Choi},\ and\ \citenamefont {Kim}}]{hong2002synchronization}%
  \BibitemOpen
  \bibfield  {author} {\bibinfo {author} {\bibfnamefont {H.}~\bibnamefont
  {Hong}}, \bibinfo {author} {\bibfnamefont {M.~Y.}\ \bibnamefont {Choi}},\
  and\ \bibinfo {author} {\bibfnamefont {B.~J.}\ \bibnamefont {Kim}},\
  }\bibfield  {title} {\bibinfo {title} {{Synchronization on small-world
  networks}},\ }\href {https://doi.org/10.1103/physreve.65.026139} {\bibfield
  {journal} {\bibinfo  {journal} {Physical Review E}\ }\textbf {\bibinfo
  {volume} {65}},\ \bibinfo {pages} {1} (\bibinfo {year} {2002})}\BibitemShut
  {NoStop}%
\bibitem [{\citenamefont {Hong}\ \emph
  {et~al.}(2007{\natexlab{a}})\citenamefont {Hong}, \citenamefont {Chaté},
  \citenamefont {Park},\ and\ \citenamefont {Tang}}]{hong2007entrainment}%
  \BibitemOpen
  \bibfield  {author} {\bibinfo {author} {\bibfnamefont {H.}~\bibnamefont
  {Hong}}, \bibinfo {author} {\bibfnamefont {H.}~\bibnamefont {Chaté}},
  \bibinfo {author} {\bibfnamefont {H.}~\bibnamefont {Park}},\ and\ \bibinfo
  {author} {\bibfnamefont {L.~H.}\ \bibnamefont {Tang}},\ }\bibfield  {title}
  {\bibinfo {title} {{Entrainment transition in populations of random frequency
  oscillators}},\ }\href {https://doi.org/10.1103/physrevlett.99.184101}
  {\bibfield  {journal} {\bibinfo  {journal} {Physical Review Letters}\
  }\textbf {\bibinfo {volume} {99}},\ \bibinfo {pages} {1} (\bibinfo {year}
  {2007}{\natexlab{a}})}\BibitemShut {NoStop}%
\bibitem [{\citenamefont {Peter}\ and\ \citenamefont
  {Pikovsky}(2018)}]{peter2018transition}%
  \BibitemOpen
  \bibfield  {author} {\bibinfo {author} {\bibfnamefont {F.}~\bibnamefont
  {Peter}}\ and\ \bibinfo {author} {\bibfnamefont {A.}~\bibnamefont
  {Pikovsky}},\ }\bibfield  {title} {\bibinfo {title} {{Transition to
  collective oscillations in finite Kuramoto ensembles}},\ }\href
  {https://doi.org/10.1103/physreve.97.032310} {\bibfield  {journal} {\bibinfo
  {journal} {Physical Review E}\ }\textbf {\bibinfo {volume} {97}},\ \bibinfo
  {pages} {032310} (\bibinfo {year} {2018})}\BibitemShut {NoStop}%
\bibitem [{\citenamefont {Brankov}\ \emph {et~al.}(2000)\citenamefont
  {Brankov}, \citenamefont {Danchev},\ and\ \citenamefont
  {Tonchev}}]{brankov2000theory}%
  \BibitemOpen
  \bibfield  {author} {\bibinfo {author} {\bibfnamefont {J.~G.}\ \bibnamefont
  {Brankov}}, \bibinfo {author} {\bibfnamefont {D.~M.}\ \bibnamefont
  {Danchev}},\ and\ \bibinfo {author} {\bibfnamefont {N.~S.}\ \bibnamefont
  {Tonchev}},\ }\href {https://doi.org/10.1142/9789812813435\_0006} {\emph
  {\bibinfo {title} {{Theory of Critical Phenomena in Finite-Size Systems}}}}\
  (\bibinfo  {publisher} {World Scientific},\ \bibinfo {year}
  {2000})\BibitemShut {NoStop}%
\bibitem [{\citenamefont {Binder}(1987)}]{binder1987finite}%
  \BibitemOpen
  \bibfield  {author} {\bibinfo {author} {\bibfnamefont {K.}~\bibnamefont
  {Binder}},\ }\bibfield  {title} {\bibinfo {title} {{Finite size effects on
  phase transitions}},\ }\href {https://doi.org/10.1080/00150198708227908}
  {\bibfield  {journal} {\bibinfo  {journal} {Ferroelectrics}\ }\textbf
  {\bibinfo {volume} {73}},\ \bibinfo {pages} {43} (\bibinfo {year}
  {1987})}\BibitemShut {NoStop}%
\bibitem [{\citenamefont {Sornette}(2006)}]{sornette2006critical}%
  \BibitemOpen
  \bibfield  {author} {\bibinfo {author} {\bibfnamefont {D.}~\bibnamefont
  {Sornette}},\ }\href {https://doi.org/10.1007/3-540-33182-4} {\emph {\bibinfo
  {title} {{Critical Phenomena in Natural Sciences, Chaos, Fractals,
  Selforganization and Disorder: Concepts and Tools}}}},\ Springer Series in
  Synergetics\ (\bibinfo {year} {2006})\BibitemShut {NoStop}%
\bibitem [{\citenamefont {Wiseman}\ and\ \citenamefont
  {Domany}(1995)}]{wiseman1995lack}%
  \BibitemOpen
  \bibfield  {author} {\bibinfo {author} {\bibfnamefont {S.}~\bibnamefont
  {Wiseman}}\ and\ \bibinfo {author} {\bibfnamefont {E.}~\bibnamefont
  {Domany}},\ }\bibfield  {title} {\bibinfo {title} {{Lack of self-averaging in
  critical disordered systems}},\ }\href
  {https://doi.org/10.1103/physreve.52.3469} {\bibfield  {journal} {\bibinfo
  {journal} {Physical Review E}\ }\textbf {\bibinfo {volume} {52}},\ \bibinfo
  {pages} {3469} (\bibinfo {year} {1995})}\BibitemShut {NoStop}%
\bibitem [{\citenamefont {Hong}\ \emph
  {et~al.}(2007{\natexlab{b}})\citenamefont {Hong}, \citenamefont {Park},\ and\
  \citenamefont {Tang}}]{hong2007finitesizescalingpre}%
  \BibitemOpen
  \bibfield  {author} {\bibinfo {author} {\bibfnamefont {H.}~\bibnamefont
  {Hong}}, \bibinfo {author} {\bibfnamefont {H.}~\bibnamefont {Park}},\ and\
  \bibinfo {author} {\bibfnamefont {L.~H.}\ \bibnamefont {Tang}},\ }\bibfield
  {title} {\bibinfo {title} {{Finite-size scaling of synchronized oscillation
  on complex networks}},\ }\href {https://doi.org/10.1103/physreve.76.066104}
  {\bibfield  {journal} {\bibinfo  {journal} {Physical Review E}\ }\textbf
  {\bibinfo {volume} {76}},\ \bibinfo {pages} {1} (\bibinfo {year}
  {2007}{\natexlab{b}})}\BibitemShut {NoStop}%
\bibitem [{\citenamefont {Skardal}\ \emph {et~al.}(2014)\citenamefont
  {Skardal}, \citenamefont {Taylor},\ and\ \citenamefont
  {Sun}}]{skardal2014optimal}%
  \BibitemOpen
  \bibfield  {author} {\bibinfo {author} {\bibfnamefont {P.~S.}\ \bibnamefont
  {Skardal}}, \bibinfo {author} {\bibfnamefont {D.}~\bibnamefont {Taylor}},\
  and\ \bibinfo {author} {\bibfnamefont {J.}~\bibnamefont {Sun}},\ }\bibfield
  {title} {\bibinfo {title} {{Optimal Synchronization of Complex Networks}},\
  }\href {https://doi.org/10.1103/physrevlett.113.144101} {\bibfield  {journal}
  {\bibinfo  {journal} {Physical Review Letters}\ }\textbf {\bibinfo {volume}
  {113}},\ \bibinfo {pages} {144101} (\bibinfo {year} {2014})}\BibitemShut
  {NoStop}%
\bibitem [{\citenamefont {Brede}(2008)}]{brede2008synchrony}%
  \BibitemOpen
  \bibfield  {author} {\bibinfo {author} {\bibfnamefont {M.}~\bibnamefont
  {Brede}},\ }\bibfield  {title} {\bibinfo {title} {{Synchrony-optimized
  networks of non-identical Kuramoto oscillators}},\ }\href
  {https://doi.org/10.1016/j.physleta.2007.11.069} {\bibfield  {journal}
  {\bibinfo  {journal} {Physics Letters, Section A: General, Atomic and Solid
  State Physics}\ }\textbf {\bibinfo {volume} {372}},\ \bibinfo {pages} {2618}
  (\bibinfo {year} {2008})}\BibitemShut {NoStop}%
\bibitem [{\citenamefont {Carareto}\ \emph {et~al.}(2009)\citenamefont
  {Carareto}, \citenamefont {Orsatti},\ and\ \citenamefont
  {Piqueira}}]{carareto2009optimized}%
  \BibitemOpen
  \bibfield  {author} {\bibinfo {author} {\bibfnamefont {R.}~\bibnamefont
  {Carareto}}, \bibinfo {author} {\bibfnamefont {F.~M.}\ \bibnamefont
  {Orsatti}},\ and\ \bibinfo {author} {\bibfnamefont {J.~R.}\ \bibnamefont
  {Piqueira}},\ }\bibfield  {title} {\bibinfo {title} {{Optimized network
  structure for full-synchronization}},\ }\href
  {https://doi.org/10.1016/j.cnsns.2008.09.032} {\bibfield  {journal} {\bibinfo
   {journal} {Communications in Nonlinear Science and Numerical Simulation}\
  }\textbf {\bibinfo {volume} {14}},\ \bibinfo {pages} {2536} (\bibinfo {year}
  {2009})}\BibitemShut {NoStop}%
\bibitem [{\citenamefont {Tang}(2011)}]{tang2011finite}%
  \BibitemOpen
  \bibfield  {author} {\bibinfo {author} {\bibfnamefont {L.~H.}\ \bibnamefont
  {Tang}},\ }\bibfield  {title} {\bibinfo {title} {{To synchronize or not to
  synchronize, that is the question: Finite-size scaling and fluctuation
  effects in the Kuramoto model}},\ }\href
  {https://doi.org/10.1088/1742-5468/2011/01/p01034} {\bibfield  {journal}
  {\bibinfo  {journal} {Journal of Statistical Mechanics: Theory and
  Experiment}\ }\textbf {\bibinfo {volume} {2011}},\ \bibinfo {pages} {P01034}
  (\bibinfo {year} {2011})},\ \Eprint {https://arxiv.org/abs/1011.5823}
  {1011.5823} \BibitemShut {NoStop}%
\bibitem [{\citenamefont {Hong}\ \emph {et~al.}(2015)\citenamefont {Hong},
  \citenamefont {Chaté}, \citenamefont {Tang},\ and\ \citenamefont
  {Park}}]{hong2015finite}%
  \BibitemOpen
  \bibfield  {author} {\bibinfo {author} {\bibfnamefont {H.}~\bibnamefont
  {Hong}}, \bibinfo {author} {\bibfnamefont {H.}~\bibnamefont {Chaté}},
  \bibinfo {author} {\bibfnamefont {L.~H.}\ \bibnamefont {Tang}},\ and\
  \bibinfo {author} {\bibfnamefont {H.}~\bibnamefont {Park}},\ }\bibfield
  {title} {\bibinfo {title} {{Finite-size scaling, dynamic fluctuations, and
  hyperscaling relation in the Kuramoto model}},\ }\href
  {https://doi.org/10.1103/physreve.92.022122} {\bibfield  {journal} {\bibinfo
  {journal} {Physical Review E}\ }\textbf {\bibinfo {volume} {92}},\ \bibinfo
  {pages} {1} (\bibinfo {year} {2015})}\BibitemShut {NoStop}%
\bibitem [{\citenamefont {Rackauckas}\ and\ \citenamefont
  {Nie}(2016)}]{rackauckas2016differential}%
  \BibitemOpen
  \bibfield  {author} {\bibinfo {author} {\bibfnamefont {C.}~\bibnamefont
  {Rackauckas}}\ and\ \bibinfo {author} {\bibfnamefont {Q.}~\bibnamefont
  {Nie}},\ }\bibfield  {title} {\bibinfo {title} {{DifferentialEquations.jl –
  A Performant and Feature-Rich Ecosystem for Solving Differential Equations in
  Julia}},\ }\href {https://doi.org/10.5334/jors.151} {\bibfield  {journal}
  {\bibinfo  {journal} {Journal of Open Research Software}\ }\textbf {\bibinfo
  {volume} {5}},\ \bibinfo {pages} {15} (\bibinfo {year} {2016})}\BibitemShut
  {NoStop}%
\bibitem [{\citenamefont {Bezanson}\ \emph {et~al.}(2017)\citenamefont
  {Bezanson}, \citenamefont {Edelman}, \citenamefont {Karpinski},\ and\
  \citenamefont {Shah}}]{bezanson2017julia}%
  \BibitemOpen
  \bibfield  {author} {\bibinfo {author} {\bibfnamefont {J.}~\bibnamefont
  {Bezanson}}, \bibinfo {author} {\bibfnamefont {A.}~\bibnamefont {Edelman}},
  \bibinfo {author} {\bibfnamefont {S.}~\bibnamefont {Karpinski}},\ and\
  \bibinfo {author} {\bibfnamefont {V.~B.}\ \bibnamefont {Shah}},\ }\bibfield
  {title} {\bibinfo {title} {{Julia: A fresh approach to numerical
  computing}},\ }\href {https://doi.org/10.1137/141000671} {\bibfield
  {journal} {\bibinfo  {journal} {SIAM Review}\ }\textbf {\bibinfo {volume}
  {59}},\ \bibinfo {pages} {65} (\bibinfo {year} {2017})}\BibitemShut {NoStop}%
\bibitem [{\citenamefont {Hunter}(2007)}]{hunter2007matplotlib}%
  \BibitemOpen
  \bibfield  {author} {\bibinfo {author} {\bibfnamefont {J.~D.}\ \bibnamefont
  {Hunter}},\ }\bibfield  {title} {\bibinfo {title} {{Matplotlib: A 2D Graphics
  Environment}},\ }\href {https://doi.org/10.1109/mcse.2007.55} {\bibfield
  {journal} {\bibinfo  {journal} {Computing in science \& engineering}\
  }\textbf {\bibinfo {volume} {9}},\ \bibinfo {pages} {90} (\bibinfo {year}
  {2007})}\BibitemShut {NoStop}%
\bibitem [{\citenamefont {Datseris}\ \emph {et~al.}(2020)\citenamefont
  {Datseris}, \citenamefont {Isensee}, \citenamefont {Pech},\ and\
  \citenamefont {Gál}}]{datseris2020drwatson}%
  \BibitemOpen
  \bibfield  {author} {\bibinfo {author} {\bibfnamefont {G.}~\bibnamefont
  {Datseris}}, \bibinfo {author} {\bibfnamefont {J.}~\bibnamefont {Isensee}},
  \bibinfo {author} {\bibfnamefont {S.}~\bibnamefont {Pech}},\ and\ \bibinfo
  {author} {\bibfnamefont {T.}~\bibnamefont {Gál}},\ }\bibfield  {title}
  {\bibinfo {title} {{DrWatson: the perfect sidekick for your scientific
  inquiries}},\ }\href {https://doi.org/10.21105/joss.02673} {\bibfield
  {journal} {\bibinfo  {journal} {Journal of Open Source Software}\ }\textbf
  {\bibinfo {volume} {5}},\ \bibinfo {pages} {2673} (\bibinfo {year}
  {2020})}\BibitemShut {NoStop}%
\bibitem [{\citenamefont {Rossi}(2022)}]{rossi2022github}%
  \BibitemOpen
  \bibfield  {author} {\bibinfo {author} {\bibfnamefont {K.}~\bibnamefont
  {Rossi}},\ }\href@noop {} {\bibinfo {title} {Repository for dynamical
  malleability code}},\ \bibinfo {howpublished}
  {\url{https://github.com/KalelR/dynamicalmalleability}} (\bibinfo {year}
  {2022})\BibitemShut {NoStop}%
\bibitem [{\citenamefont {Medvedev}(2014)}]{medvedev2014small}%
  \BibitemOpen
  \bibfield  {author} {\bibinfo {author} {\bibfnamefont {G.~S.}\ \bibnamefont
  {Medvedev}},\ }\bibfield  {title} {\bibinfo {title} {{Small-world networks of
  Kuramoto oscillators}},\ }\href {https://doi.org/10.1016/j.physd.2013.09.008}
  {\bibfield  {journal} {\bibinfo  {journal} {Physica D: Nonlinear Phenomena}\
  }\textbf {\bibinfo {volume} {266}},\ \bibinfo {pages} {13} (\bibinfo {year}
  {2014})}\BibitemShut {NoStop}%
\bibitem [{\citenamefont {Carlson}\ \emph {et~al.}(2011)\citenamefont
  {Carlson}, \citenamefont {Kim},\ and\ \citenamefont
  {Motter}}]{carlson2011sample}%
  \BibitemOpen
  \bibfield  {author} {\bibinfo {author} {\bibfnamefont {N.}~\bibnamefont
  {Carlson}}, \bibinfo {author} {\bibfnamefont {D.-H.}\ \bibnamefont {Kim}},\
  and\ \bibinfo {author} {\bibfnamefont {A.~E.}\ \bibnamefont {Motter}},\
  }\bibfield  {title} {\bibinfo {title} {{Sample-to-sample fluctuations in
  real-network ensembles}},\ }\href {https://doi.org/10.1063/1.3602223}
  {\bibfield  {journal} {\bibinfo  {journal} {Chaos: An Interdisciplinary
  Journal of Nonlinear Science}\ }\textbf {\bibinfo {volume} {21}},\ \bibinfo
  {pages} {025105} (\bibinfo {year} {2011})},\ \Eprint
  {https://arxiv.org/abs/1111.6118} {1111.6118} \BibitemShut {NoStop}%
\bibitem [{\citenamefont {Hong}\ \emph {et~al.}(2006)\citenamefont {Hong},
  \citenamefont {Park},\ and\ \citenamefont {Tang}}]{hong2006anomalous}%
  \BibitemOpen
  \bibfield  {author} {\bibinfo {author} {\bibfnamefont {H.}~\bibnamefont
  {Hong}}, \bibinfo {author} {\bibfnamefont {H.}~\bibnamefont {Park}},\ and\
  \bibinfo {author} {\bibfnamefont {L.-H.}\ \bibnamefont {Tang}},\ }\bibfield
  {title} {\bibinfo {title} {{Anomalous Binder Cumulant and Lack of
  Self-Averageness in Systems with Quenched Disorder}},\ }\href@noop {}
  {\bibfield  {journal} {\bibinfo  {journal} {Journal of the Korean Physical
  Society}\ }\textbf {\bibinfo {volume} {49}},\ \bibinfo {pages} {L1885}
  (\bibinfo {year} {2006})}\BibitemShut {NoStop}%
\bibitem [{\citenamefont {Strogatz}\ and\ \citenamefont
  {Mirollo}(1988)}]{strogatz1988phase}%
  \BibitemOpen
  \bibfield  {author} {\bibinfo {author} {\bibfnamefont {S.~H.}\ \bibnamefont
  {Strogatz}}\ and\ \bibinfo {author} {\bibfnamefont {R.~E.}\ \bibnamefont
  {Mirollo}},\ }\bibfield  {title} {\bibinfo {title} {{Phase-locking and
  critical phenomena in lattices of coupled nonlinear oscillators with random
  intrinsic frequencies}},\ }\href
  {https://doi.org/10.1016/0167-2789(88)90074-7} {\bibfield  {journal}
  {\bibinfo  {journal} {Physica D: Nonlinear Phenomena}\ }\textbf {\bibinfo
  {volume} {31}},\ \bibinfo {pages} {143} (\bibinfo {year} {1988})}\BibitemShut
  {NoStop}%
\bibitem [{\citenamefont {Tilles}\ \emph {et~al.}(2011)\citenamefont {Tilles},
  \citenamefont {Ferreira},\ and\ \citenamefont
  {Cerdeira}}]{tilles2011multistable}%
  \BibitemOpen
  \bibfield  {author} {\bibinfo {author} {\bibfnamefont {P.~F.~C.}\
  \bibnamefont {Tilles}}, \bibinfo {author} {\bibfnamefont {F.~F.}\
  \bibnamefont {Ferreira}},\ and\ \bibinfo {author} {\bibfnamefont {H.~A.}\
  \bibnamefont {Cerdeira}},\ }\bibfield  {title} {\bibinfo {title}
  {{Multistable behavior above synchronization in a locally coupled Kuramoto
  model}},\ }\href {https://doi.org/10.1103/physreve.83.066206} {\bibfield
  {journal} {\bibinfo  {journal} {Physical Review E}\ }\textbf {\bibinfo
  {volume} {83}},\ \bibinfo {pages} {066206} (\bibinfo {year}
  {2011})}\BibitemShut {NoStop}%
\bibitem [{\citenamefont {Budzinski}\ \emph {et~al.}(2020)\citenamefont
  {Budzinski}, \citenamefont {Rossi}, \citenamefont {Boaretto}, \citenamefont
  {Prado},\ and\ \citenamefont {Lopes}}]{budzinski2020synchronization}%
  \BibitemOpen
  \bibfield  {author} {\bibinfo {author} {\bibfnamefont {R.~C.}\ \bibnamefont
  {Budzinski}}, \bibinfo {author} {\bibfnamefont {K.~L.}\ \bibnamefont
  {Rossi}}, \bibinfo {author} {\bibfnamefont {B.~R.~R.}\ \bibnamefont
  {Boaretto}}, \bibinfo {author} {\bibfnamefont {T.~L.}\ \bibnamefont
  {Prado}},\ and\ \bibinfo {author} {\bibfnamefont {S.~R.}\ \bibnamefont
  {Lopes}},\ }\bibfield  {title} {\bibinfo {title} {{Synchronization
  malleability in neural networks under a distance-dependent coupling}},\
  }\href {https://doi.org/10.1103/physrevresearch.2.043309} {\bibfield
  {journal} {\bibinfo  {journal} {Physical Review Research}\ }\textbf {\bibinfo
  {volume} {2}},\ \bibinfo {pages} {043309} (\bibinfo {year}
  {2020})}\BibitemShut {NoStop}%
\bibitem [{\citenamefont {Newman}\ and\ \citenamefont
  {Watts}(1999)}]{newman1999scaling}%
  \BibitemOpen
  \bibfield  {author} {\bibinfo {author} {\bibfnamefont {M.~E.}\ \bibnamefont
  {Newman}}\ and\ \bibinfo {author} {\bibfnamefont {D.~J.}\ \bibnamefont
  {Watts}},\ }\bibfield  {title} {\bibinfo {title} {{Scaling and percolation in
  the small-world network model.}},\ }\href
  {https://doi.org/10.1103/physreve.60.7332} {\bibfield  {journal} {\bibinfo
  {journal} {Physical Review E}\ }\textbf {\bibinfo {volume} {60}},\ \bibinfo
  {pages} {7332} (\bibinfo {year} {1999})}\BibitemShut {NoStop}%
\bibitem [{\citenamefont {Villegas}\ \emph {et~al.}(2014)\citenamefont
  {Villegas}, \citenamefont {Moretti},\ and\ \citenamefont
  {Muñoz}}]{villegas2014frustrated}%
  \BibitemOpen
  \bibfield  {author} {\bibinfo {author} {\bibfnamefont {P.}~\bibnamefont
  {Villegas}}, \bibinfo {author} {\bibfnamefont {P.}~\bibnamefont {Moretti}},\
  and\ \bibinfo {author} {\bibfnamefont {M.~A.}\ \bibnamefont {Muñoz}},\
  }\bibfield  {title} {\bibinfo {title} {{Frustrated hierarchical
  synchronization and emergent complexity in the human connectome network}},\
  }\href {https://doi.org/10.1038/srep05990} {\bibfield  {journal} {\bibinfo
  {journal} {Scientific Reports}\ }\textbf {\bibinfo {volume} {4}},\ \bibinfo
  {pages} {5990} (\bibinfo {year} {2014})},\ \Eprint
  {https://arxiv.org/abs/1402.5289} {1402.5289} \BibitemShut {NoStop}%
\bibitem [{\citenamefont {Tononi}\ \emph {et~al.}(1994)\citenamefont {Tononi},
  \citenamefont {Sporns},\ and\ \citenamefont {Edelman}}]{tononi1994a}%
  \BibitemOpen
  \bibfield  {author} {\bibinfo {author} {\bibfnamefont {G.}~\bibnamefont
  {Tononi}}, \bibinfo {author} {\bibfnamefont {O.}~\bibnamefont {Sporns}},\
  and\ \bibinfo {author} {\bibfnamefont {G.~M.}\ \bibnamefont {Edelman}},\
  }\bibfield  {title} {\bibinfo {title} {{A measure for brain complexity:
  relating functional segregation and integration in the nervous system.}},\
  }\href {https://doi.org/10.1073/pnas.91.11.5033} {\bibfield  {journal}
  {\bibinfo  {journal} {Proceedings of the National Academy of Sciences}\
  }\textbf {\bibinfo {volume} {91}},\ \bibinfo {pages} {5033} (\bibinfo {year}
  {1994})}\BibitemShut {NoStop}%
\bibitem [{\citenamefont {Deco}\ \emph {et~al.}(2015)\citenamefont {Deco},
  \citenamefont {Tononi}, \citenamefont {Boly},\ and\ \citenamefont
  {Kringelbach}}]{deco2015rethinking}%
  \BibitemOpen
  \bibfield  {author} {\bibinfo {author} {\bibfnamefont {G.}~\bibnamefont
  {Deco}}, \bibinfo {author} {\bibfnamefont {G.}~\bibnamefont {Tononi}},
  \bibinfo {author} {\bibfnamefont {M.}~\bibnamefont {Boly}},\ and\ \bibinfo
  {author} {\bibfnamefont {M.~L.}\ \bibnamefont {Kringelbach}},\ }\bibfield
  {title} {\bibinfo {title} {{Rethinking segregation and integration:
  contributions of whole-brain modelling}},\ }\href
  {https://doi.org/10.1038/nrn3963} {\bibfield  {journal} {\bibinfo  {journal}
  {Nature Reviews Neuroscience}\ }\textbf {\bibinfo {volume} {16}},\ \bibinfo
  {pages} {430} (\bibinfo {year} {2015})}\BibitemShut {NoStop}%
\bibitem [{\citenamefont {Tononi}\ and\ \citenamefont
  {Edelman}(1998)}]{tononi1998consciousness}%
  \BibitemOpen
  \bibfield  {author} {\bibinfo {author} {\bibfnamefont {G.}~\bibnamefont
  {Tononi}}\ and\ \bibinfo {author} {\bibfnamefont {G.~M.}\ \bibnamefont
  {Edelman}},\ }\bibfield  {title} {\bibinfo {title} {{Consciousness and
  Complexity}},\ }\href {https://doi.org/10.1126/science.282.5395.1846}
  {\bibfield  {journal} {\bibinfo  {journal} {Science}\ }\textbf {\bibinfo
  {volume} {282}},\ \bibinfo {pages} {1846} (\bibinfo {year}
  {1998})}\BibitemShut {NoStop}%
\bibitem [{\citenamefont {Fingelkurts}\ and\ \citenamefont
  {Fingelkurts}(2006)}]{fingelkurts2006timing}%
  \BibitemOpen
  \bibfield  {author} {\bibinfo {author} {\bibfnamefont {A.~A.}\ \bibnamefont
  {Fingelkurts}}\ and\ \bibinfo {author} {\bibfnamefont {A.~A.}\ \bibnamefont
  {Fingelkurts}},\ }\bibfield  {title} {\bibinfo {title} {{Timing in cognition
  and EEG brain dynamics: discreteness versus continuity.}},\ }\href
  {https://doi.org/10.1007/s10339-006-0035-0} {\bibfield  {journal} {\bibinfo
  {journal} {Cognitive Processing}\ }\textbf {\bibinfo {volume} {7}},\ \bibinfo
  {pages} {135} (\bibinfo {year} {2006})}\BibitemShut {NoStop}%
\bibitem [{\citenamefont {Li}\ \emph {et~al.}(2009)\citenamefont {Li},
  \citenamefont {Poo},\ and\ \citenamefont {Dan}}]{li2009burst}%
  \BibitemOpen
  \bibfield  {author} {\bibinfo {author} {\bibfnamefont {C.-y.~T.}\
  \bibnamefont {Li}}, \bibinfo {author} {\bibfnamefont {M.-m.}\ \bibnamefont
  {Poo}},\ and\ \bibinfo {author} {\bibfnamefont {Y.}~\bibnamefont {Dan}},\
  }\bibfield  {title} {\bibinfo {title} {{Burst Spiking of a Single Cortical
  Neuron Modifies Global Brain State}},\ }\href
  {https://doi.org/10.1126/science.1169957} {\bibfield  {journal} {\bibinfo
  {journal} {Science}\ }\textbf {\bibinfo {volume} {324}},\ \bibinfo {pages}
  {643} (\bibinfo {year} {2009})}\BibitemShut {NoStop}%
\bibitem [{\citenamefont {Muller}\ \emph {et~al.}(2021)\citenamefont {Muller},
  \citenamefont {Min{\'a}{\v{c}}},\ and\ \citenamefont
  {Nguyen}}]{muller2021algebraic}%
  \BibitemOpen
  \bibfield  {author} {\bibinfo {author} {\bibfnamefont {L.}~\bibnamefont
  {Muller}}, \bibinfo {author} {\bibfnamefont {J.}~\bibnamefont
  {Min{\'a}{\v{c}}}},\ and\ \bibinfo {author} {\bibfnamefont {T.~T.}\
  \bibnamefont {Nguyen}},\ }\bibfield  {title} {\bibinfo {title} {{Algebraic
  approach to the Kuramoto model}},\ }\href
  {https://doi.org/10.1103/physreve.104.l022201} {\bibfield  {journal}
  {\bibinfo  {journal} {Physical Review E}\ }\textbf {\bibinfo {volume}
  {104}},\ \bibinfo {pages} {L022201} (\bibinfo {year} {2021})}\BibitemShut
  {NoStop}%
\bibitem [{\citenamefont {Budzinski}\ \emph {et~al.}(2022)\citenamefont
  {Budzinski}, \citenamefont {Nguyen}, \citenamefont {Do{\`a}n},
  \citenamefont {Min{\'a}{\v{c}}}, \citenamefont {Sejnowski},\ and\
  \citenamefont {Muller}}]{budzinski2022geometry}%
  \BibitemOpen
  \bibfield  {author} {\bibinfo {author} {\bibfnamefont {R.~C.}\ \bibnamefont
  {Budzinski}}, \bibinfo {author} {\bibfnamefont {T.~T.}\ \bibnamefont
  {Nguyen}}, \bibinfo {author} {\bibfnamefont {J.}~\bibnamefont
  {Do{\`a}n}}, \bibinfo {author} {\bibfnamefont {J.}~\bibnamefont
  {Min{\'a}{\v{c}}}}, \bibinfo {author} {\bibfnamefont {T.~J.}\ \bibnamefont
  {Sejnowski}},\ and\ \bibinfo {author} {\bibfnamefont {L.~E.}\ \bibnamefont
  {Muller}},\ }\bibfield  {title} {\bibinfo {title} {{Geometry unites
  synchrony, chimeras, and waves in nonlinear oscillator networks}},\ }\href
  {https://doi.org/10.1063/5.0078791} {\bibfield  {journal} {\bibinfo
  {journal} {Chaos}\ }\textbf {\bibinfo {volume} {32}},\ \bibinfo {pages}
  {031104} (\bibinfo {year} {2022})},\ \Eprint
  {https://arxiv.org/abs/2111.02560} {2111.02560} \BibitemShut {NoStop}%
\end{thebibliography}
%

\end{document}